\documentclass[11pt]{extarticle}
\usepackage[utf8]{inputenc}
\usepackage[T1,T2A]{fontenc}
\usepackage{geometry}
\usepackage{amsmath,amsfonts,amssymb,amscd,euscript}
\usepackage{amsthm}
\usepackage[dvipsnames]{xcolor}
\usepackage[english]{babel}
\usepackage{graphicx} 
\usepackage{psfrag} 
\usepackage{float}
\usepackage{subfig}
\usepackage{setspace}
\usepackage{bm}
\usepackage{mathabx}
\usepackage{url}
\usepackage{multirow,tabularx}

\usepackage[ruled,vlined]{algorithm2e}

\graphicspath{{Figures/}}

\numberwithin{equation}{section}

\title{A Bayesian Approach to Modelling Biological Pattern Formation with Limited Data}

\author{Alexey Kazarnikov$^1$, Robert Scheichl$^1$, Heikki Haario$^2$ and Anna Marciniak-Czochra$^{1}$}

\date{\footnotesize
${}^1$ Institute of Applied Mathematics and Interdisciplinary Center for Scientific Computing (IWR), Heidelberg University, Im Neuenheimer Feld 205, 69120 Heidelberg, Germany\\[0.5ex]
${}^2$ Department of Computational Engineering, LUT University, Yliopistonkatu 34, 53850 Lappeenranta, Finland\\[1ex]
}

\begin{document}

\def\Re{\mathop{\rm Re}\,}
\def\Im{\mathop{\rm Im}\,}
\def\dom{\mathop{\rm dom}\,}
\def\dist{\mathop{\rm dist}}
\def\grad{\mathop{\rm grad}}
\newcommand{\phan}{\hspace*{0cm}}
\newcommand{\comment}{}
\newcommand\No{\char"9D}

\newtheorem{proposition}{Proposition}
\newtheorem{corollary}{Corollary}
\newtheorem{definition}{Definition}
\newtheorem{remark}{Remark}
\newtheorem{Algorithm}{Algorithm}

\makeatletter
\g@addto@macro\@floatboxreset\centering
\makeatother





\maketitle

\begin{abstract}
Understanding the mechanisms of pattern formation is one of the fundamental questions in developmental biology. Since the classic work of Alan Turing, reaction-diffusion systems have been the dominant modelling approach. Alternative models combine the dynamics of diffusing molecular signals with tissue mechanics or intracellular feedback. 
Model validation is complicated as in many experimental situations only the limiting, stationary regime of the pattern formation process can be observed, without any knowledge of the transient behaviour or the initial state. To overcome this problem, the initial state of the model can be randomised. However, then the fixed values of the model parameters correspond to a family of patterns rather than a fixed stationary solution, and standard estimation approaches such as least squares are not appropriate.  Instead, statistical characteristics of the  patterns should be compared, which is computationally expensive and difficult given the limited amount of data usually available in practical applications.
To deal with this problem, we extend a recently developed statistical approach for parameter identification by pattern data, the Correlation Integral Likelihood (CIL) method.  We introduce modifications of the CIL approach that allow to increase the accuracy of the identification process without resizing the data set, in particular a range of alternative distance measures as well as the synthetic likelihood approach. The proposed modifications are tested using different classes of pattern formation models and severely limited data sets. The numerical solvers for all considered equations are based on highly scalable, parallel, GPU-based implementations with efficient time stepping schemes.
\end{abstract}

\section{Introduction}

Mechanisms of self-organised pattern formation in developmental biology have been in the focus of experimental and theoretical research for several decades \cite{Lewis2008,Meinhardt2006,murray:1,tomlin2007}.  A number of different morphologies have been modelled mathematically \cite{kondo:1995,Nakamasu2009,Torii2012}. All those models, developed under different biological hypotheses, have common mathematical features and are mostly based on a limited range of mathematical paradigms of pattern generation. The classical Turing approach relies on a chemical pre-pattern created by hypothetical morphogens \cite{turing:1}. Nonlinear interactions, in conjunction with different diffusivity, may yield a bifurcation  (diffusion-driven instability), which leads to destabilisation of a spatially homogeneous equilibrium and to the emergence of patterns. Stable patterns can also result from the coupling of a single diffusive morphogen to a non-diffusive
subsystem, exhibiting a hysteresis effect  resulting in {\it far-from-equilibrium} patterns \cite{harting2017,halatek2018,kothe2020}. An alternative mechanism of pattern formation has been emphasised by mathematical models linking the dynamics of diffusing molecular signals with tissue mechanics \cite{mercker2013,mercker2015, goriely2017, murray1984, murray1984a, oster1985}. 

Various theoretical models hint at how complex intracellular pathways coupled with cell-to-cell communication through diffusing molecules or mechanical signal transduction can lead to symmetry breaking and spatially heterogeneous patterns. However, qualitative and quantitative differences between the different models remain unexplored. It is unknown whether they produce similar stable patterns, or whether the patterns are similarly sensitive to changes in biophysical parameters { or initial conditions}. Answers to these mathematical questions will help designing experiments to validate the nature of distinct pattern formation steps during development and regeneration.

In order to identify a model, it should be verified against experimental observations. However, in many experimental situations only a limiting, stationary regime of the pattern formation process is observable and data on transient behaviour or the initial state are missing. 
As a starting point for model identification, it is already interesting and challenging enough to use synthetic data produced by the models, to then estimate parameter values solely from { the time-limiting stationary patterns. Since
different dynamics may lead to similar stationary patterns \cite{oster1988, economou2014, sharpe2019, woolley2021}, 
it is conceivable that in certain cases stationary patterns are fundamentally insufficient for inference. A successful methodology should then allow identifying alternative models, which could be further investigated to find possible differences in the dynamics, for example, in response to certain perturbations. Thus, such an approach would help to constrain the experimental search for an explanation of the underlying mechanisms on a molecular level.} 

Another important difficulty associated with numerical identification of pattern formation models results from the pattern selection problem.
Due to instabilities in the underlying dynamics, even small changes in the initial condition may lead to distinct limiting stationary patterns, which seriously affects parameter identification. { Such a phenomenon is observed in models of {\it de novo} pattern formation, such as the classical Turing models or mechanical-chemical models, and has been studied by many authors in the context of model robustness \cite{maini2012}. 
Coexistence of distinct patterns can also occur due to hysteresis effects, 
leading to the existence of different stable branches of quasi-stationary solutions, see \cite{Veerman2021} for a brief comparison of {\it close-to-equilibrium} and {\it far-from-equilibrium} patterns. Altogether, pattern selection is a complex phenomenon that must be carefully considered during computational model identification.} 

To overcome the instability problem, the initial state in the model can be randomised. Accordingly, fixed values of the model parameters correspond to a family of patterns rather than to a fixed stationary solution. Standard approaches to compare pattern data directly with model outputs, e.g., in the least squares sense, are meaningless. Therefore, many existing approaches for parameter identification employ either known initial values or transient data \cite{kramer:1,garvie:1,Campillo-Funollet2019}. Additional difficulties arise from numerical approximation errors that may affect the transient behaviour and the final limiting state. Consequently, a change in numerical resolution may lead to an entirely different pattern, even with identical initial conditions.
Furthermore, data are often available only in normalised form, providing information about the pattern shape but not quantitative concentration levels.


These challenges have motivated a new approach to model calibration, based on a statistical comparison of patterns on the basis of a stochastic cost function that quantifies the correlation between point clouds, the so-called \emph{Correlation Integral Likelihood (CIL)} \cite{kazarnikov2020}. The CIL approach provides a likelihood that is normally distributed. Thus, many available optimisation algorithms can be employed to compute the maximum likelihood estimate. Standard Bayesian techniques are available to estimate the posterior distributions and thus to quantify the uncertainty due to measurement noise and due to a limited amount of available data. In \cite{kazarnikov2020}, the concept was applied to three reaction-diffusion models exhibiting Turing patterns. However, the proposed approach requires a lot of data to construct the likelihood. It is estimated that at least 50 data patterns, each containing about five ‘wavelengths’, are necessary \cite{kazarnikov2020}. The prediction accuracy significantly drops as the size of the data set is decreased. In the limiting case of 50 patterns, values ‘within’ and ‘outside’ the posterior distribution can be distinguished even by the naked eye.

Recently, also neural network approaches  
have been applied to Turing models using stationary pattern data. In \cite{zhu2022}, convolutional neural networks are employed to learn the data-to-parameter map for SIR-type rumour propagation models, using a sufficiently large training set and a suitable neural network architecture. A comparison of the approach in \cite{zhu2022} with the CIL approach in \cite{kazarnikov2020} suggests that CIL requires significantly more data than neural networks. Additionally, a loss of accuracy of the CIL method was observed in \cite{zhu2022}. However, this latter observation depends on the choice of data, the specific model and on its implementation. In \cite{schnorr2021}, a range of machine learning approaches were compared on an example of the Gierer-Meinhardt model. The identification process was significantly improved by introducing an invariant pattern representation based on resistance distance histograms, which made it possible to represent the spatial structure of patterns and to compute distances between pattern data effectively. Once trained with a sufficiently large data set, the methods in \cite{zhu2022,schnorr2021} allowed to identify model parameters accurately from a single pattern. { However, the pattern data in \cite{zhu2022} stems from a significantly larger domain containing patterns with about 
10-15 ‘wavelengths’, while the pattern data in \cite{kazarnikov2020} contains only about 4-5 ‘wavelengths’.} 
  
In the current work, the applicability of the CIL approach is significantly expanded. Addressing the criticism raised in  \cite{zhu2022,schnorr2021}, we test the method on a diverse set of pattern formation models. First, we consider cases with very limited data sets, resembling practical applications. Then, we propose modifications increasing the accuracy of the identification without resizing the data set. A key innovation is the introduction of several distance measures. In addition to the $L^2$-norm employed in  \cite{kazarnikov2020}, we evaluate distances between patterns using  $L^\infty$-, $H^1$- and $W^{1,\infty}$-norms. In general, any number of additional ‘features’ could be added to this list, { including problem-specific measures not based on distances}. The introduction of multiple norms results in a significant contraction of the posterior distributions of the parameters, and thus a significant reduction in the residual uncertainty. We conduct parameter identification using pattern data sets that are severely limited, down to the case of a single pattern with less than 5 ‘wavelengths’. This is made possible by our second key innovation: the combination of the CIL likelihood with the \emph{synthetic likelihood} sampling approach \cite{wood:1,price:1}. The method is tested on different classes of pattern formation models, including the Turing reaction-diffusion systems, mechano-chemical models, as well as reaction-diffusion-ODE systems exhibiting \emph{far-from-equilibrium patterns}. 

{ For all considered equations, parallel GPU-based implementations of the numerical solvers with improved time stepping schemes are introduced that significantly reduce the computational cost compared to the implementation of CIL in \cite{kazarnikov2020} (and certainly to its implementation in \cite{zhu2022}). A batched parallel treatment of multiple model trajectories on GPUs allows for an efficient online inference of model parameters for all models without costly pre-training.} Finally, the approach is also compared to the neural network approach and the benefits of each method are discussed. 

The text is organised as follows. Main concepts of the CIL algorithm are recalled in Section~\ref{sec:para_ident}. Section \ref{sec:modif} introduces the multi-feature idea and the reformulation of the CIL approach as a synthetic likelihood method to successfully tackle the limited data case. In Section~\ref{sec:models}, the different models are presented, followed by numerical experiments in Section~\ref{sec:numerics}. The description of the spatial and temporal discretisation schemes are relegated to Appendix~\ref{sec:discretise}. The paper concludes with a discussion of the results and a comparison to other recent approaches in the literature in Section~\ref{sec:discussion}. 

\section{Parameter identification by pattern data}
\label{sec:para_ident}

We consider the problem of model parameter identification by pattern data. Only the limiting, stationary regime of the pattern formation process is assumed to be observed experimentally, without any detailed knowledge of the transient behaviour or the initial data. 
%
The main difficulty is in choosing a suitable cost function, which quantifies the goodness of fit of a model to the available data for given parameter values. We solve this problem by a statistical approach, which allows defining a stochastic cost function that makes it possible to measure a distance between experimental pattern data and the model output. This technique was developed in \cite{haario:1} and \cite{springer2021} for classical, chaotic dynamical systems and extended to reaction-diffusion systems in \cite{kazarnikov2020}. Once defined, the CIL cost function can be minimised by available algorithms of stochastic optimisation, resulting in point estimates for model parameters. One can continue with Bayesian methods to statistically quantify the identification accuracy. These key stages of the approach are explained in detail below.


\subsection{Correlation Integral Likelihood: A cost function for parameter identification}

Let us consider an abstract pattern formation model depending on a vector of control parameters $\bm{\theta} \in \mathbb{R}^p$. We assume that the output of the model are patterns $\bm{s}(\bm{\theta})  { \in \mathbb{R}^d}$. These patterns naturally change due to the variation of model parameters, but additionally change even for fixed model parameters due to instabilities in the nonlinear pattern formation process, e.g., due to small variations in the initial conditions. In what follows, we implicitly assume that the patterns being produced are presented as finite-dimensional approximations of the limiting, function-valued stationary numerical solution of a nonlinear differential equation model (see below).

Our goal is to construct an empirical statistical likelihood for the parameter $\bm{\theta}$ from a set of patterns to be used for parameter identification. Various methods based on summary statistics have been developed for 'intractable' situations  where explicit statistical likelihoods are not available \cite{wood:1,price:1,sisson2020}.  Here, we employ an approach based on generalisations of the Central Limit Theorem (CLT). According to the CLT, averages of bounded observations can be approximated by Gaussian distributions. Averaging, however, leads to a loss of information. Instead, the empirical distribution of data can be used. { The empirical cumulative distribution function  (eCDF) of independent and identically distributed (i.i.d.) scalar random values is computed by 
binning the data, giving eCDF values $0< p_i <1$, $i=1,...,M$. The probability of a random value to be below or above $p_i$ is then Bernoulli distributed (with expected value $p_i$ and variance $p_i(1-p_i) )$, and the average of such values tends to a Gaussian distribution by the CLT. The resulting vector, which contains the eCDF values at all bin values, tends to a multidimensional Gaussian distribution, with  the dimension equal to the number of bins.
More generally, the basic form of Donsker's theorem states that the eCDF 
of i.i.d.\ scalar random variables asymptotically tends to a Gaussian process  \cite{donsker1951,donsker1952}. 
%
 In our setting the data is not i.i.d., but the Gaussianity is established by theorems of the U-statistics \cite{borovkova:1,neumeyer:1} that allow weakly dependent observations.}
 
 Nevertheless, as the pattern data is high dimensional, a mapping to scalar values is needed before creating the eCDF vectors. This mapping may again lead to a loss of information and has to be selected with care. In \cite{haario:1,springer2019,kazarnikov2020}, the Correlation Integral Likelihood (CIL), a modification of the correlation fractal dimension concept, was successfully used for this purpose. The general idea of this procedure is recalled below, while in Section \ref{sec:multi}, the approach is further extended and refined to allow for several scalar-valued mappings in order to more extensively employ the pattern characteristics of the original data. In Section \ref{sec:synthetic}, we will further advance the approach to allow for parameter identification using the CIL also in the case of very small data sets.
 
 

 
Let $\bm{s}_\text{data} = \{\bm{s}_i : i=1,\ldots,N_\text{set}\}$ be a set of $N_\text{set}$ patterns in $\mathbb{R}^d$, $d \in \mathbb{N}$, that are assumed to all come from the same underlying pattern formation process, in particular from the same model with fixed model parameter $\bm{\theta}_0 \in \mathbb{R}^p$, $p \in \mathbb{N}$.  For any two subsets $\bm{s}$ and $\widetilde{\bm{s}}$ of $\bm{s}_\text{data}$ containing $N$ and $\widetilde{N}$ patterns, respectively, and for a fixed 'radius' $R>0$, we define
\begin{equation}
	C(R,\bm{s},\widetilde{\bm{s}},\|\cdot\|)=\frac{1}{N \times \widetilde{N}}\sum_{i=1}^{N} \sum_{j=1}^{\widetilde{N}} \#(\|\bm{s}_{i}-\widetilde{\bm{s}}_{j}\|<R),
	\label{eq:eCDF}
\end{equation}
where $\|\cdot\|$ is a suitable norm. Now, to analyse the data set $\bm{s}_\text{data}$ we first subdivide it into $n_\text{ens}$ subsets, denoted $\bm{s}^k$, $k=1,\ldots,n_\text{ens}$, each consisting of $N$ samples of patterns, such that $N_\text{set} = n_\text{ens} \times N$. Next, we choose $R_{0}>0$ such that  $\|\bm{s}^k_i-\bm{s}^l_j\| < R_0$, for all $k,l = 1,\ldots,n_\text{ens}$, $i,j=1,\ldots,N$, $k\ne l$, and define $R_m := g(R_0,m)$ with a function $g$ that is strictly decreasing in the second argument describing the decay rate of the radii. For a given integer $M$, the 
vector $\bm{y}^{k,l}\in\mathbb{R}^M$ defined componentwise by 
\begin{equation}
	y_{m}^{k,l}=C(R_m,\bm{s}^k,\bm{s}^l,\|\cdot\|),\quad m=1,2,\dots,M,
	\label{eq:CIL}
\end{equation}
then defines a realisation of the so-called \emph{correlation integral vector} $\bm{y}_0$, the eCDF vector of the distances of two sets of $N$ patterns $\bm{s}^k$ and $\bm{s}^l$, both coming from the same model with identical parameters $\bm{\theta}_0$ and evaluated at the bin values $R_m$.

In the numerical experiments, the largest and smallest radii, $R_0$ and $R_M$, are determined according to the range of all the distances in (\ref{eq:eCDF}). Function $g(R_0,m)$ is either the simple power law relationship $g(R_0,m) = R_0 b^{-m}$ with $b>1$ such that $R_M/R_0 = b^{-M}$, or the linear relationship $g(R_0,m) = R_0 - m h$ with $ h = (R_0-R_M)/M$. 

%
To define the \emph{Correlation Integral Likelihood (CIL)} of an arbitrary parameter $\bm{\theta}\in \mathbb{R}^p$, we first numerically estimate the mean $\bm{\mu}_0\in\mathbb{R}^M$ and covariance $\bm{\Sigma}_0\in\mathbb{R}^{M\times M}$ of the correlation integral vector $\bm{y}_0$, using $\left({n_\text{ens} \atop 2}\right)$ realisations defined in \eqref{eq:CIL} for all pairs of subsets of the $N_\text{set}$ given patterns with $k \ne l$. As stated above,  under suitable conditions it has been shown rigorously that $\bm{y}_0$ follows a multivariate Gaussian distribution $N(\bm{\mu}_0,\bm{\Sigma}_0)$. { The main assumption is that the data is only weakly dependent, which may be difficult to verify theoretically. So in practice we test numerically for Gaussianity of the ensemble of vectors defined in \eqref{eq:CIL} using the $\chi^2$-test (or a scalar normality test for each of the components of the vector).} 
Note that the construction of the mean and covariance in \eqref{cil:CF} can be performed off-line, without any further model runs, using the training data. It is summarised in Algorithm \ref{alg:CIL1}. { However, when the amount of data is limited the amount of patterns per subset is also severely limited, so that it is difficult to guarantee a sufficiently large number of pairs for the stable estimation of $\bm{\mu}_0$ and $\bm{\Sigma}_0$. This naturally increases the variability of $\bm{y}_0$, which in turn affects the accuracy of the parameter identification.} 
\begin{algorithm}[t]
\SetKwInOut{Input}{Input}
\SetKwInOut{Output}{Output}
\SetKwInOut{Data}{Data}
\SetKwInOut{Result}{Result}
\Data{A set $\bm{s}_\text{data}$ of discretised patterns $\bm{s}_1,\ldots,\bm{s}_{N_\text{set}}$, with unknown model parameter $\bm{\theta}_0$}
\Result{Parameters $\bm{\mu}_0$ and $\bm{\Sigma}_0$ for the CIL cost function $f(\bm{\theta})$ in \eqref{cil:CF}}   	
\Begin
{
1. Divide $\bm{s}_\text{data}$ into $n_\text{ens}$ subsets $\bm{s}^k$, $k=1,\ldots,n_\text{ens}$, each consisting of $N$ patterns

2. \For{$k,l=1,\ldots,n_\mathrm{ens}, k\neq l$}
 {
 	* Compute the $N^2$ distances $\|\bm{s}^k_{i}-\bm{s}^l_{j}\|$ between all patterns in $\bm{s}^k$ and $\bm{s}^l$ in norm $\|\cdot\|$ 
 	
 	* Compute a sample of the correlation integral vector $\bm{y}^{k,l}$ as defined in \eqref{eq:CIL}
 }

 3. Using the samples $\bm{y}^{k,l}$ computed in Step 2, estimate mean $\bm{\mu}_0$ and covariance matrix~$\bm{\Sigma}_0$\!\!\!\!\!\linebreak of the correlation integral vector $\bm{y}_0$ of the training data.
 
 }
 \caption{Construction of the cost function in the basic CIL approach}
 \label{alg:CIL1}
 \end{algorithm}

Now, to quantify the 'distance' between patterns in the training set $\bm{s}_\text{data}$ and patterns produced by the model with parameter $\bm{\theta} \in \mathbb{R}^p$, we define a \emph{generalised correlation integral vector} $\bm{y}(\bm{\theta})$ as follows: a set $\bm{s}(\bm{\theta})$ of $N$ patterns $\bm{s}_i(\bm{\theta})$, $i=1,\ldots,N$, is computed from $N$ runs of the model with parameter~$\bm{\theta}$; this is then compared with a randomly selected  subset $\bm{s}^k$ of the original training data to give
\begin{equation}
	y_{m}(\bm{\theta})=C(R_m,\bm{s}^k,\bm{s}(\bm{\theta}),\|\cdot\|),\quad m=1,2,\dots,M.
	\label{eq:CIL_theta}
\end{equation}
Alternative definitions of $\bm{y}(\bm{\theta})$ are possible, e.g., by averaging the distances of $\bm{s}(\bm{\theta})$ to all the subsets $\bm{s}^k$, $k=1,\ldots,n_\text{ens}$. Under the assumption that the CIL of $\bm{y}_0$ is Gaussian $N(\bm{\mu_0},\bm{\Sigma_0})$, the cost function for parameter estimation can now be defined as the negative log-likelihood function
\begin{equation}
	f(\bm{\theta}) = \big(\bm{y}(\bm{\theta})-\bm{\mu}_{0} \big)^\top \bm{\Sigma}_{0}^{-1} \big(\bm{y}(\bm{\theta})-\bm{\mu}_{0}\big).
	\label{cil:CF}
\end{equation}
A more detailed description of the whole procedure can be found in \cite{kazarnikov2020}. 

To find the minimum of the stochastic cost function $f(\bm{\theta})$ in \eqref{cil:CF} different methods of stochastic optimisation 
can be employed. The previous work \cite{kazarnikov2020} used Differential Evolution (DE), a stochastic, population-based 
algorithm \cite{storn:1}. The procedure starts from a randomly distributed population of parameter vectors on a bounded set in parameter space. The population evolves iteratively 
by following prescribed rules, and eventually converges to a local minimum of the objective function.

\subsection{Uncertainty quantification}

After obtaining the maximum a posteriori (MAP) point of the parameter distribution using the DE algorithm, an ensemble of samples from the posterior distribution can be constructed by classical Markov Chain Monte Carlo (MCMC) sampling \cite{Robert2004}. One can generate a 'chain' of parameter values $\bm{\theta}^1$, $\bm{\theta}^2$, $\ldots$, $\bm{\theta}^n$, whose empirical distribution  approaches the posterior distribution. 
Starting from the computed MAP point, we use DRAM \cite{MCMC:2}, an adaptive Metropolis-Hastings algorithm \cite{MCMC:1} with random walk proposals, where the covariance matrix of the proposal distribution is updated during sampling. The proposal covariance can be initialised using the final population of the DE algorithm. The generalised correlation integral vector $\bm{y}(\bm{\theta})$ is stochastic, as the initial conditions for the patterns $\bm{s}_i(\bm{\theta})$, $i=1,\ldots,N$, are random and also the subset $\bm{s}^k$ that the patterns are compared against are drawn at random. Thus, the MCMC method can be interpreted as a pseudomarginal sampler. 

\section{Modifications of the parameter identification approach}
\label{sec:modif}

The CIL approach is based on the assumption that the number of training patterns $N_\text{set}$ is sufficiently large to produce reliable estimates for the parameters of the underlying Gaussian distribution. In many practical applications, the required amount of data might be prohibitively large. This restriction can be weakened by using resampling methods, such as bootstrapping, suggested and successfully employed in \cite{kazarnikov2020}. Bootstrapping worked in a stable way to estimate the parameters for $N_\text{set} \geq 50$, even though the accuracy of the parameter identification decreased with the size of the training set.
Smaller data sets did not provide reliable estimates for the CIL parameters.

In the current paper, we suggest two innovations that allow further reducing the size of data sets. The first modification significantly improves the accuracy of the parameter identification without increasing the amount of training data. The second modification allows extending the method to limited data sets, where the small number of available patterns prevents the usage of resampling methods. These modifications are discussed in detail below.

\subsection{Multi-feature Correlation Integral Likelihood}
\label{sec:multi}

The accuracy of parameter identification achieved by the CIL approach decreases with the amount of available data. Here, we establish a significant improvement of the accuracy of the method without increasing the training set. The general idea is to use several mappings from pattern data to scalar values in the CIL scheme (discussed in the previous section), thus highlighting different features in the patterns and reducing the loss of information. In particular, we use different norms in the definition of the generalised correlation integral vectors in \eqref{eq:CIL} and \eqref{eq:CIL_theta}.
Since these vectors are assumed to follow Gaussian distributions, any number $K$ of different features can be simply concatenated into one larger generalised correlation integral vector $\bm{y}_0 \in \mathbb{R}^{MK}$, which still follows a multivariate Gaussian distribution $N(\bm{\mu_0},\bm{\Sigma_0})$ on $\mathbb{R}^{MK}$. For general $\bm{\theta}$, the stochastic cost function $f(\bm{\theta})$ in \eqref{cil:CF} also remains unchanged. The additional computational overhead due to applying multiple scalar mappings to the simulated patterns is negligible.

To leverage different data features, we use different distances of the spatial patterns. This is motivated by applicability of different functional setting in nonlinear stability analysis of  reaction-diffusion equations \cite{KowallMarciniakMuennich2022}.  Specifically, we use discrete equivalents of the following norms\\[-3ex]
%
\hspace{-0.5cm}\begin{minipage}[t]{0.47\textwidth}
\begin{eqnarray}
 \|u(\bm{x})\|_{L^2} &=&\bigg({\int_{\Omega}|u(\bm{x})|^{2} \text{d}x}\bigg)^{1/2}, \quad\label{norm:L2} \\[0.5ex]
 \|u(\bm{x})\|_{L^\infty} &= &\sup_{\bm{x}\in\Omega}|u(\bm{x})|,\quad \label{norm:Linf} \\[0.5ex]
 \vvvert u(\bm{x})\vvvert_{W^{1,2}} &=& \sum_{|\alpha|\leq1} \|D^{\alpha}(u(\bm{x}))\|_{L^{2}}, \quad\label{norm:W12prime}
\end{eqnarray}
\end{minipage}
\begin{minipage}[t]{0.53\textwidth}
\begin{eqnarray}
 \|u(\bm{x})\|_{W^{1,2}} &=& \bigg(\sum\nolimits_{|\alpha|\leq1}\|D^{\alpha}(u(\bm{x}))\|_{L^{2}}^2\bigg)^{1/2}\!\!, \quad\label{norm:W12} \\[0.5ex]
 \|u(\bm{x})\|_{W^{1,\infty}} &=& \max_{|\alpha|\leq1} \|D^{\alpha}(u(\bm{x}))\|_{L^{\infty}}, \quad\label{norm:W1inf} \\[0.5ex]
 \vvvert u(\bm{x}) \vvvert_{W^{1,\infty}} &=& \sum_{|\alpha|\leq1}\|D^{\alpha}(u(\bm{x}))\|_{L^{\infty}}, \quad\label{norm:W1infprime}
\end{eqnarray}
\end{minipage}\\[1ex]
to quantify the difference $u(\bm{x})$ between two spatial patterns over all points $\bm{x}\in\Omega$ for some $\Omega \subset\mathbb{R}^1$ or $\Omega \subset\mathbb{R}^2$. In case of the Sobolev norms in \eqref{norm:W12prime}-\eqref{norm:W1infprime}, the derivatives  are approximated by finite differences in our numerical experiments.
Other alternative norms or mappings are conceivable. For example, it might be useful to compute distances between patterns on suitably chosen subdomains.


The procedure to compute the new \emph{Multi-feature Correlation Integral Likelihood (MCIL)} is summarised in Algorithm \ref{alg:CIL2}. 

Further improvement of the MCIL-based parameter estimates can be achieved with bootstrapping, similar as in the case of the standard CIL approach in \cite{kazarnikov2020}. It applies especially to cases with limited data. However, also for moderately-sized or large data sets (as considered in \cite{zhu2022}, where $N_{set} = 3000$), bootstrapping significantly reduces the variability of the MCIL parameter estimates. Hence, the new approach solves the problem of variability of the CIL criticised in \cite{zhu2022}. For completeness, the MCIL approach with bootstrapping is presented in Algorithm \ref{alg:CIL2bst} in Appendix \ref{sec:algos}.

\begin{algorithm}[t]
	\SetKwInOut{Input}{Input}\SetKwInOut{Output}{Output}
	\SetKwInOut{Data}{Data}
	\SetKwInOut{Result}{Result}
	\Data{A set $\bm{s}_\text{data}$ of discretised patterns $\bm{s}_1,\ldots,\bm{s}_{N_\text{set}}$, with unknown model parameter $\bm{\theta}_0$}
	\Input{Family of norms $\|\cdot\|_\alpha$, $\alpha=1,\ldots,N_\text{dist}$}
	\Result{Parameters $\bm{\mu}_0$ and $\bm{\Sigma}_0$ for the cost function $f(\bm{\theta})$ in \eqref{cil:CF} in the case of the MCIL}
	\Begin
	{
		1. Divide $\bm{s}_\text{data}$ into $n_\text{ens}$ subsets $\bm{s}^k$, $k=1,\ldots,n_\text{ens}$, each consisting of $N$ patterns
		
		2. \For{$k,l=1,\ldots,n_\mathrm{ens}, k\neq l$}
		{
		2.1 Initialise the correlation integral vector $\bm{y}^{k,l}$ for $\bm{s}^k$ and $\bm{s}^l$ to be the empty vector

		2.2 \For{$\alpha=1,\ldots,N_\mathrm{dist}$}
		{
			* Using $\|\cdot\|_\alpha$ compute the distances $\|\bm{s}^k_{i}-\bm{s}^l_{j}\|_\alpha$ between all patterns in $\bm{s}^k$ and $\bm{s}^l$\!\!
			
			* Using $\|\cdot\|_\alpha$ compute the part $\bm{y}^{k,l}_\alpha$ of the correlation integral vector via \eqref{eq:CIL}
			
			* Concatenate the current vector $\bm{y}^{k,l}$ and $\bm{y}^{k,l}_\alpha$
		}
    }
	
    3. Using the samples $\bm{y}^{k,l}$ computed in Step 2, estimate mean $\bm{\mu}_0$ and covariance matrix~$\bm{\Sigma}_0$\!\!\!\!\!\linebreak of the (multi-feature) correlation integral vector $\bm{y}_0$ of the training data
	}
	\caption{Construction of the cost function in the multi-feature CIL approach}
	\label{alg:CIL2}
\end{algorithm}

\subsection{Synthetic Correlation Integral Likelihood}
\label{sec:synthetic}

Next, we consider a situation with a limited training data set, such that the resampling techniques are not sufficient, i.e., when the set $\bm{s}_\text{data} = \{\bm{s}_i:i=1,\ldots,N_\text{set}\}$ produced with the unknown parameter $\bm{\theta}_0$ is restricted to $N_\text{set} = O(10)$ or even $O(1)$ patterns. In that case, we propose a different approach: To define a cost function $f(\bm{\theta})$ for any given parameter value $\bm{\theta}$, we approximate the CIL (or the MCIL) at $\bm{\theta}$ using a sufficiently large set $\{\bm{s}_i(\bm{\theta}) : i=1,\ldots,N_\text{syn}\}$ of patterns, computed from $N_\text{syn}$ model runs, and then compare it to the training data $\bm{s}_\text{data}$. This idea is a special case of the \emph{synthetic likelihood method}, introduced in \cite{wood:1} and further discussed in, e.g., \cite{price:1}.

In particular, in contrast to the standard CIL algorithm~\ref{alg:CIL1} and to the multi-feature CIL algorithm~\ref{alg:CIL2}, here we subdivide the (synthetic) patterns $\bm{s}_i(\bm{\theta})$, $i=1,\ldots,N_\text{syn}$, at the parameter value $\bm{\theta}$ into subsets $\bm{s}^k_{\bm{\theta}}$, $k=1,\ldots,n_\text{ens}$, with $N$ patterns each, such that $N_\text{syn} = n_\text{ens} \times N$ and $N > N_\text{set}$. 
The subsets are then used to define realisations  $\widetilde{\bm{y}}^{k,l}_{\bm{\theta}}$ of a (synthetic) correlation integral vector at $\bm{\theta}$, similar to the one in \eqref{eq:CIL}. However, we need to modify the definition slightly. In the \emph{synthetic correlation integral likelihood (SCIL)} we want to reliably estimate the 'distance' of the data set $\bm{s}_\text{data}$ to the synthetic pattern data, even when its dimension $N_\text{set}$ is very small (in the worst case just $N_\text{set}=1$). Thus, to match the likelihood construction and the evaluation of it with $\bm{s}_\text{data}$, each subset $\bm{s}^k_{\bm{\theta}}$, $k=1,\ldots,n_\text{ens}$, is further partitioned into two subsets $\bm{s}^{k,1}_{\bm{\theta}}$ and $\bm{s}^{k,2}_{\bm{\theta}}$ with $N_\text{set}$ and $\widetilde{N} := N - N_\text{set}$ patterns, respectively. We define
\begin{equation}
	\widetilde{y}_{\bm{\theta},m}^{k,l}=C(R_m,\bm{s}^{k,1},\bm{s}^{l,2},\|\cdot\|),\quad m=1,2,\dots,M.
	\label{eq:CIL_synth_kl}
\end{equation}
Note that in contrast to \eqref{eq:CIL} the two sets of patterns are of different dimensions, but since the patterns are distinct even when $k=l$, all $n_\text{ens}^2$ combinations of $\bm{s}^{k,1}$ and $\bm{s}^{l,2}$ are taken into account. To reliably estimate the bin values of the eCDF of the 'distances' of the patterns in $\bm{s}^{k,1}$ and $\bm{s}^{l,2}$ in \eqref{eq:CIL_synth_kl}, $N$ needs to be chosen sufficiently large, i.e., such that $N_\text{set} \times \widetilde{N}$ is sufficiently large. The Gaussianity of the distribution of the vectors $\widetilde{\bm{y}}^{k,l}_{\bm{\theta}}$ can be verified again numerically using $\chi^2$-test.

The radii $R_m$, $m=1,\ldots,M$, are chosen adaptively by first determining the maximum (resp. minimum) distances $R_0$ (resp. $R_M$) of any two patterns and then fitting the power law decay rate $b>1$ such that $R_M/R_0 = b^{-M}$ or the linear progression $R_m = R_0 - m h$ with $h = (R_0-R_M)/M$. Note, however, that since the range of distances does in general vary with $\bm{\theta}$, the radii and thus also the bin values for the eCDF of the correlation integral vectors will vary with $\bm{\theta}$.

Finally, we can estimate the mean $\bm{\mu}_{\bm{\theta}}$ and covariance matrix $\bm{\Sigma}_{\bm{\theta}}$ of the (synthetic) correlation integral vector $\widetilde{\bm{y}}_{\bm{\theta}}$ to obtain the objective function for the SCIL at $\bm{\theta}$: 
\begin{equation}
	f(\bm{\theta}) = \big(\widetilde{\bm{y}}(\bm{\theta})-\bm{\mu}_{\bm{\theta}} \big)^\top \bm{\Sigma}_{\bm{\theta}}^{-1} \big(\widetilde{\bm{y}}(\bm{\theta})-\bm{\mu}_{\bm{\theta}}\big),\vspace{-1ex}
	\label{cil:CF_synth}
\end{equation}
where 
\begin{equation}
	\widetilde{y}_{m}(\bm{\theta})=C(R_m,\bm{s}_\text{data},\bm{s}^{{ k_0},2},\|\cdot\|),\quad m=1,2,\dots,M,
	\label{eq:CIL_synth}
\end{equation}
is the generalised correlation integral vector. Here, $\bm{s}_\text{data}$ is the entire training data set at the unknown parameter $\bm{\theta}_0$ and { $k_0$ is a subset index chosen uniformly at random from $\{1,\ldots,n_{ens}\}$.} 

\begin{algorithm}[t]
	\SetKwInOut{Input}{Input}\SetKwInOut{Output}{Output}
	\SetKwInOut{Data}{Data}
	\SetKwInOut{Result}{Result}
	\Data{A set $\bm{s}_\text{data}$ of discretised patterns $\bm{s}_1,\ldots,\bm{s}_{N_\text{set}}$, with unknown model parameter $\bm{\theta}_0$}
	\Input{$\|\cdot\|_\alpha  \ \ldots$ \ a family of norms with $\alpha=1,\ldots,N_\text{dist}$}
	\Input{$\bm{\theta} \quad\ \ \ \ldots$ \ the parameter value where the cost function should be evaluated}
	\Input{$N_\text{syn} \ \ \ldots$ \ the number of model-generated (synthetic) patterns at $\bm{\theta}$}
	\Output{$f(\bm{\theta}) \ \ \ldots$ \ the value of the  cost function at $\bm{\theta}$}	
	\Begin
	{
		1. Simulate $N_\text{syn}$ patterns $\bm{s}_i(\bm{\theta})$, $i=1,\ldots,N_\text{syn}$, with parameter $\bm{\theta}$
		
		2. Divide these patterns into $2 \times n_\text{ens}$ subsets $\bm{s}_{\bm{\theta}}^{k,1}$ and $\bm{s}_{\bm{\theta}}^{k,2}$, $k=1,\ldots,n_\text{ens}$, with $N_\text{set}$ and $\widetilde{N}$ patterns, respectively
		
		3. \For{$k,l=1,\ldots,n_\mathrm{ens}$}
		{
		3.1 Initialise the correlation integral vector $\bm{y}_{\bm{\theta}}^{k,l}$ for $\bm{s}_{\bm{\theta}}^{k,1}$ and $\bm{s}_{\bm{\theta}}^{l,2}$ to be the empty vector\!\!\! 

		3.2 \For{$\alpha=1,\ldots,N_\mathrm{dist}$}
		    {
			* Compute the distances $\|\bm{s}^{k,1}_{\bm{\theta},i}-\bm{s}^{l,2}_{\bm{\theta},j}\|_\alpha$ between all patterns in $\bm{s}_{\bm{\theta}}^{k,1}$ and $\bm{s}_{\bm{\theta}}^{l,2}$ w.~$\|\cdot\|_\alpha$\!\!\!\!
			
			* Compute part of correlation integral vector $\bm{y}^{k,l}_{\bm{\theta},\alpha}$ as defined in \eqref{eq:CIL_synth_kl} using $\|\cdot\|_\alpha$
			
			* Concatenate the current vector $\bm{y}_{\bm{\theta}}^{k,l}$ and $\bm{y}^{k,l}_{\bm{\theta},\alpha}$
		    }
        }		
        4. Using the samples $\bm{y}_{\bm{\theta}}^{k,l}$ computed in Step 3, estimate mean $\bm{\mu}_{\bm{\theta}}$ and covariance matrix~$\bm{\Sigma}_{\bm{\theta}}$\!\!\!\!\!\linebreak of the (multi-feature) correlation integral vector $\bm{y}_{\bm{\theta}}$ of the synthetic data at $\bm{\theta}$

		5. Randomly select a subset $\bm{s}_{\bm{\theta}}^{k,2}$ of $\widetilde{N}$ patterns from
		$\bm{s}_i(\bm{\theta})$, $i=1,\ldots,N_\text{syn}$

		6. Using the computed estimates $\bm{\mu}_{\bm{\theta}}$ and $\bm{\Sigma}_{\bm{\theta}}$, the set $\bm{s}_{\bm{\theta}}^{k,2}$ and $\widetilde{\bm{y}}(\theta)$ as defined in~\eqref{eq:CIL_synth}, evaluate $f(\bm{\theta})$,  the cost function at $\bm{\theta}$ from \eqref{cil:CF_synth}
	}
	\caption{Construction and evaluation of  $f(\bm{\theta})$ for SCIL at a single parameter value $\bm{\theta}$}
	\label{alg:CIL3}
\end{algorithm}

Multi-feature identification, as discussed in the previous section, can naturally be applied in the current context as well. The complete procedure is summarised in Algorithm \ref{alg:CIL3}. In addition, bootstrapping may be applied to lower the required amount of model simulations and thus  to significantly decrease the cost without any loss of accuracy (see Algorithm \ref{alg:CIL3bst} in Appendix \ref{sec:algos}).
To obtain stable estimates for the mean vector $\bm{\mu}_{\bm{\theta}}$ and covariance matrix $\bm{\Sigma}_{\bm{\theta}}$ of the CIL (or MCIL) at each evaluation of the cost function $f(\bm{\theta})$, a significantly larger number of model runs are required than for CIL or MCIL. More precisely, the overall cost of the synthetic likelihood approach requires $O(n_\text{ens})$ times more model runs per evaluation of $f(\bm{\theta})$ than the original CIL or the MCIL approach above.
However, this additional overhead is acceptable in the significantly more challenging, limited data case, and due to our efficient parallel GPU implementation the total effort remains reasonable.

\section{Test cases: Models for biological pattern formation}
\label{sec:models}

We consider three classes of models exhibiting diverse pattern formation mechanisms: classical reaction-diffusion systems, one-dimensional mechano-chemical models, and reaction-diffusion-ODE systems.  { All the reaction-diffusion systems are considered on the unit square $\Omega=(0,1)^2$, while the two other models are posed on the one-dimensional interval $\Omega=(0,1)$. All equations are defined in non-dimensionalised form and subject to homogeneous Neumann (zero-flux) boundary conditions.}

\subsection{Reaction-diffusion systems}

Classical Turing-type models of pattern formation are given in the form of two-component reaction-diffusion systems,
\begin{equation}
	v_t = \nu_1\Delta v + f(v,w), \qquad
	w_t = \nu_2\Delta w + g(v,w), \label{eq:rdgen}
\end{equation}
where $v=v(\bm{x},t)$ and $w=w(\bm{x},t)$ describe concentrations of signalling molecules, $\nu_{1},\,\nu_{2}>0$ are  diffusion coefficients, and the non-linear functions $f(v,w)$ and $g(v,w)$ represent local interactions. 

As test cases, we choose three classical reaction-diffusion systems: the FitzHugh-Nagumo model  \cite{fitzhugh:1,nagumo:1}, the Gierer-Meinhardt activator-inhibitor system \cite{gierer:1} and the Brusselator 
model~\cite{prigogine:1}. The FitzHugh-Nagumo model provides a two-component reduction of the Hodgkin-Huxley nerve impulse propagation model and is given by the following reaction terms:
\begin{equation}
	f(v,w)=\varepsilon(w-\alpha v), \qquad
	g(v,w)=-v+\mu w-w^{3}, \label{introduction:0}
\end{equation}
where $v(\bm{x},t)$ is the recovery variable, $w(\bm{x},t)$ is the membrane potential, and $\mu\in \mathbb{R}$, $\alpha\geq0$ and $\varepsilon>0$ are control parameters.
The Gierer-Meinhardt system is a prototype short-range activation and long-range inhibition model. The reaction terms read: 
\begin{equation}
	f(v,w)= - \mu_v v + \frac{v^2}{w}, \qquad
	g(v,w)= - \mu_w w + v^2, \label{gmSystem}
\end{equation}
where $v(\bm{x},t)$ and $w(\bm{x},t)$ denote activator and inhibitor concentrations and $\mu_v,\mu_w>0$ the rescaled decay rates. 
The Brusselator is a theoretical model of autocatalytic chemical reactions, described by the following kinetics:
\begin{equation}
	f(v,w)= A - (B+1)v + v^2{w}, \qquad
	g(v,w)= Bv - v^2{w},
	\label{bzSystem}
\end{equation}
where 
$A$ and $B$ describe a constant reactant supply and a constant inverse reaction rate, respectively.

We focus on the model dynamics around a spatially homogeneous steady state that is losing stability due to the diffusion-driven instability (DDI), 
where
\[
f(v_0,w_0)=g(v_0,w_0)=0.
\]
It holds: $(v_0,w_0)=(0,0)$ for model \eqref{introduction:0}, $(v_0,w_0)=(\frac{\mu_w}{\mu_v},\frac{\mu_w}{\mu_v^2})$ for system \eqref{gmSystem} and $(v_0,w_0)=(A,\frac{B}{A})$ for equations \eqref{bzSystem}. 
Each system is considered for parameter values satisfying DDI conditions, which are computed from linearised stability analysis.

{
Thus, the initial conditions for the reaction-diffusion models in all the numerical simulations below are chosen as small uniform random noise perturbations of the steady state $(v_0,w_0)$, i.e. 
\[
v(\bm{x},0)=v_0 + U(0,\delta),
\,
w(\bm{x},0)=w_0 + U(0,\delta),
\ \text{with} \ \delta=10^{-2}.
\]
}

\subsection{Mechano-chemical models}

Mechano-chemical models may exhibit a spontaneous pattern formation that is governed by a positive feedback loop between the tissue curvature and a morphogen production. 
 We restrict ourselves to one-dimensional mechano-chemical models, which are simplified versions of more realistic equations studied in \cite{mercker2013,mercker2015}.  A biological tissue is represented by a deforming incompressible one-dimensional surface $u(x,t)$, where $x \in (0,L)$ and $t \geq 0$. By $\kappa(u)$ we denote the local tissue curvature:
\begin{equation}
\kappa(u) = \dfrac{\partial}{\partial x}\left( \dfrac{u_x(x,t)}{\sqrt{1+u_x(x,t)^2}} \right)
=
\dfrac{u_{xx}(x,t)}{\left(1+u_x(x,t)^2\right)^{\frac{3}{2}}}\;.
\label{eq:curvature}
\end{equation}
Evolution of the tissue is described by a 4th order partial differential equation coupled to a system of reaction-diffusion equations describing dynamics of biochemical morphogens $\bm{\phi}$:  
\begin{eqnarray}
\tau u_t(x,t) & =&  -L
\left[
\dfrac{1}{\sqrt{1+u_x(x,t)^2}}\dfrac{\partial}{\partial x}\left(\dfrac{\kappa_x(u)-\bar{\kappa}_x(\bm{\phi})}{\sqrt{1+u_x(x,t)^2}}\right)
-
\left(\kappa(u)-\bar{\kappa}(\bm{\phi})\right)\kappa^2(u) + \lambda\kappa(u)
\right],
\label{eq:tissue} \\
 (\phi_i(x,t))_t &=& \dfrac{D_i}{\sqrt{1+u_x(x,t)^2}}\dfrac{\partial}{\partial x}\left(\dfrac{(\phi_i(x,t))_x}{\sqrt{1+u_x(x,t)^2}}\right)+F_i(u,\bm{\phi}), \quad i=1,\ldots,n,\label{eq:morphogen}
\end{eqnarray}
where $\bar{\kappa}(\bm{\phi}) : \mathbb{R}^n\rightarrow\mathbb{R}$ denotes the locally preferred curvature, $L$ is the length of the spatial domain, $\tau$ is a relaxation parameter, $D_i$ are diffusion coefficients, $F_i(u,\bm{\phi})$ are non-linear reaction terms and $\lambda$ is the Lagrange multiplier.
 The model is supplemented by the tissue incompressibility requirement,
 \begin{eqnarray}
 S(u)
&=&
\int\limits_0^1{\sqrt{1+u_x(x,t)^2}}dx 
= 
\int\limits_0^1{\sqrt{1+u_x(x,0)^2}}dx
=
\text{const},
\label{eq:constraint}
\end{eqnarray}

We focus on the model with $n=1$, i.e.,
\begin{equation}
    \bm{\phi} \equiv \phi(x,t),
    \quad
    \bm{F}(u,\bm{\phi}) = -\alpha\phi(x,t) + f(\kappa(u)),
    \quad
    \bar{\kappa}(\bm{\phi}) = -\beta\phi(x,t),
    \label{eq:mcone}
\end{equation}
where $\alpha>0$ is the morphogen degradation rate. The morphogen production rate depends on $\kappa(u)$ and is given by  $f(\kappa(u))=\eta\dfrac{\Delta \kappa(u) 
}{1+ \Delta\kappa(u)}$ 
, where  $\Delta \kappa(u) := \max((\kappa_0-\kappa(u)),0)$ is the deviation from the initial curvature $\kappa_0 = \kappa(u(x,0))$ and $\eta > 0$ is the limit production rate.
%
The system is closed with homogeneous Neumann boundary conditions.
{
As initial data, 
we choose
\[
\phi(x,0)\equiv 0, \quad u(x,0)=1 + U(0,\delta)
\ \ \text{with} \ \  \delta=10^{-2}.
\]
}

\subsection{Reaction-diffusion-ODE models}

The third class of the pattern formation models considered in this paper is given by systems of reaction-diffusion equations coupled to space-dependent ordinary differential equation.   While such models may exhibit DDI, it turns out that all branching stationary solutions (Turing patterns) are unstable \cite{Marciniak-CzochraKarchSuzuki2017, cygan2022} and {\it far-from-equilibrium} patterns may emerge \cite{harting2017,kothe2020}. 
We consider an example model consisting of two equations\vspace{-1ex}
\begin{equation}
		u_t = -u - uw + m_1 \dfrac{u^2}{1+ku^2}\;, \qquad
		w_t = D \dfrac{\partial^2}{\partial x^2}w - m_3w - uw + m_2 \dfrac{u^2}{1+ku^2}\;,
	\label{eq:rdode}
\end{equation}
with zero-flux boundary conditions for $w$.
{Here, $u=u(x,t)$ and $w=w(x,t)$ denote concentrations of the non-diffusive and the diffusive components, respectively,
$D>0$ is the diffusion coefficient, and $m_1,m_2,m_3>0$ and $k>0$ are reaction parameters.}
{The initial data for the non-diffusive component $u(x,t)$ are created by smoothing uniform random noise using a repetitive application of a moving averaging, while the diffusive component $w(x,t)$ is started from zero.}

\section{Numerical experiments}
\label{sec:numerics}

To illustrate the effect of the proposed modifications on the accuracy of the CIL approach, we conduct numerical experiments, which are described now. 

All experiments are run for two types of pattern data. In the first case, we assume that actual concentration values of the patterns $\bm{s}(\bm{\theta})$ are known, thus we assume that the direct model output is available for observation. In the second case, we consider scaled pattern data, where the absolute values of the concentrations are removed from the patterns by a min-max normalisation:
\[
\bm{s}(\bm{\theta}) = (s_1(x),\ldots,s_n(x)) \rightarrow \tilde{\bm{s}}(\bm{\theta})  = (\tilde{s}_1(x),\ldots,\tilde{s}_n(x)), \quad
\tilde{s}_i(x) = \dfrac{s_i(x) - s_i^{min}}{s_i^{max} - s_i^{min}},
\]
where $s_{i}^{min}$ and $s_{i}^{max}$ are minimum and maximum values of the corresponding functions on the respective spatial domain $\Omega$. This situation can be considered to be representative of experimental patterns being given by greyscale pictures only. As a result of such a transformation, naturally some information is lost and the accuracy of the parameter identification is reduced. Nevertheless, we consider this case to be of great practical importance in situations where actual concentration values cannot be measured and only the shapes of the final patterns are known.

\begin{table}[t]
\small
\begin{tabular}{|c|c|c|c|c|c|c|c|}
\hline
     \textbf{Parameter name} & \textbf{FHN} & \textbf{GM} & \textbf{BZ} & \textbf{MC} & \textbf{RD}  \\\hline
     Time integration interval & $[0,100]$ & $[0,100]$ & $[0,100]$ & $[0,100]$ & $[0,50]$ \\
     Spatial resolution & $64 \times 64$ nodes & $64 \times 64$ nodes & $64 \times 64$ nodes & $64$ nodes & $64$ nodes \\
     $N_\text{set}$ (training set size) & 50 & 50 & 50 & 200 & 200 \\
     $M$ (CIL vector dimension) & 13 & 13 & 13 & 13 & 13 \\
     Radii decay law & power law & power law & power law & linear & linear \\
     $\bm{\theta}$ (control parameters) & $(\mu,\varepsilon)$ & $(\mu_v,\mu_w)$ & $(A,B)$ & $(D,\alpha)$ & $(m_1,m_2)$ \\
     $\bm{\theta}_0$ (''true'' parameter values) & $(1,10)$ & $(0.5,1)$ & $(4.5,6.96)$ & $(0.75,1)$ & $(1.44,4.1)$ \\
      & $\nu_1=0.05$ & $\nu_1 = 0.00025$ & $\nu_1=0.0016$ & $\beta=10$ & $D=1$ \\
      & $\nu_2=0.00028$ & $\nu_2 = 0.01$ & $\nu_2 = 0.0132$ & $\delta=10$ & $m_3 = 4.1$ \\
      & $\alpha=1$ & & & $L=10$ & $k = 0.01$ \\
      & & & & $\tau = 200$ & \\
     \hline
\end{tabular}
\caption{Parameter values used in the numerical experiments. The abbreviations FHN, GM, BZ, MC, and RD denote respectively the  FitzHugh-Nagumo model, the Gierer-Meinhardt system, the Brusselator 
system, the mechanochemical model 
and the reaction-diffusion-ODE model.}
\label{table:exp1}
\end{table} 

{
In this section, when referring to a pattern $\bm{s}(\bm{\theta})$, we always mean the finite-dimensional vector of (spatial) grid values, obtained as the limiting stationary numerical solution of the properly discretised nonlinear pattern formation model. The spatial discretisations of the models and the time integration schemes are discussed in detail in Appendix \ref{sec:discretise}. The time integration interval, the spatial resolution and all other relevant parameters for the numerical experiments are given in Table \ref{table:exp1}. 

As stated above, initial conditions for the reaction-diffusion models and for the mechano-chemical models are small uniform perturbations of a spatially homogeneous steady state.  We verified numerically that the distribution and the size of the random perturbations have no impact on the pattern formation process, at least for small values of $\delta$, which are reasonable for {\it de novo} pattern formation. To this end, we fixed model parameters and used different combinations of Gaussian and uniform random noise. Moreover, we used several variants of  Mat\'ern fields \cite{lindgren2011,rasmussen2005} as  they produce initial conditions varying from stochastic noise to smooth spatially inhomogeneous perturbations. Again, no impact was observed. 
In contrast, the choice of initial data does have a pronounced effect for the reaction-diffusion-ODE model and hence becomes a part of the parameter identification process.
}

In the following experiments, we first demonstrate the increase in accuracy achieved by employing the multi-feature CIL (MCIL) approach for sufficiently large data sets, before studying the performance of the synthetic likelihood idea (SCIL) in the limited data case. To optimise computational times, in both cases the bootstrapping procedure is used (Algorithms \ref{alg:CIL2bst} and \ref{alg:CIL3bst} in Appendix~\ref{sec:algos}). Throughout the experiments, we use $25$ (resp.~$100$) model runs for each evaluation of the cost function $f(\bm{\theta})$ for the reaction-diffusion models (resp.~for the mechano-chemical and reaction-diffusion-ODE models) in the case of Algorithm \ref{alg:CIL2bst}, and $1000$ model runs in the case of Algorithm \ref{alg:CIL3bst}. 
To make this amount of model runs feasible during parameter identification or posterior sampling, it is absolutely crucial to optimise the performance of the numerical solver for the forward problems.

\subsection{Multi-GPU implementation and numerical efficiency}

A key ingredient for the success of the CIL, and in particular the SCIL approach, is the efficient parallel implementation of the numerical solvers on modern Graphical Processing Units (GPUs), using the 
Nvidia CUDA (Compute Unified Device Architecture) computing platform to execute computations on massively parallel Nvidia GPU devices. However, an optimal choice and implementation of the numerical schemes has also been essential.

As we are dealing with non-linear problems, the number of operations per time step can differ significantly. Due to the internal structure of GPU hardware, it is very inefficient to run single simulations independently. Instead, computationally expensive vector operations are synchronised between trajectories and executed by the GPU device in a batched manner, such as right hand side evaluations or linear system solves. { For batched execution, it is also crucial (or at least advantageous in practice) to use uniform spatial grids.} Non-expensive scalar operations are executed on the host (CPU). All solvers have a built-in feature to scale also to multiple GPU devices if available. 

For reaction-diffusion systems, a major efficiency gain over the results in~\cite{kazarnikov2020} has been achieved by using the numerically stable, second-order, explicit time-stepping approach ROCK2 \cite{abdulle2001} instead of the explicit Euler scheme (for details see Section \ref{sec:ROCK2}). It is ideally suited for efficient batch computations on GPU platforms.  The method allows for large time steps, which results in a much smaller number of right-hand side evaluations, { yet it does not require any linear solves like implicit time stepping methods and thus scales significantly better on GPU platforms (see the comparison in Section \ref{sec:ROCK2})}. We integrate reaction-diffusion models and the reaction-diffusion-ODE system on the intervals $[0,100]$ and $[0,50]$, respectively. For the considered range of parameter domain, this is sufficient to reach convergence to the steady-state pattern.

\begin{figure}[t]
	\centering
	\includegraphics[width=0.99\textwidth]{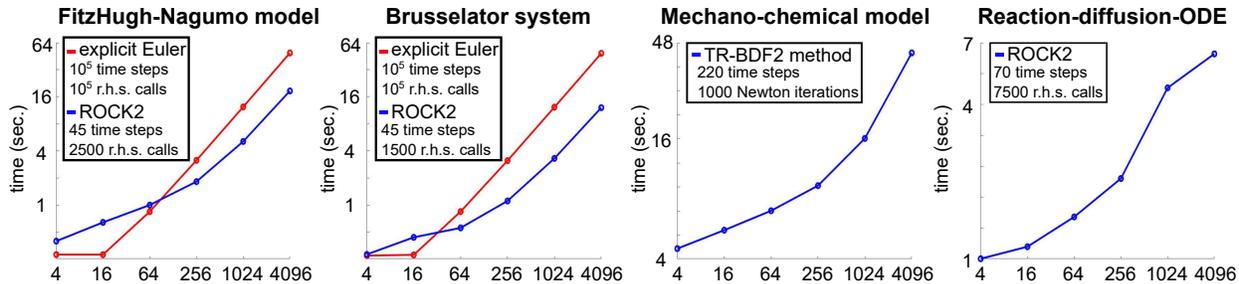}
	\caption{Execution times of the GPU-based implementations of the numerical solvers with respect to the number $n_{\text{sim}}$ of parallel simulations (logarithmic scales). 
	The experiment was performed on a single Nvidia GTX 980 graphics card.}
	\label{Fig:6}
\end{figure}

For the mechano-chemical system, which involves a fourth-order, nonlinear PDE, implicit time-stepping is essential to avoid severe time step restrictions. We use an adaptive, second-order Runge-Kutta method that combines a trapezoidal and BDF2 substep \cite{bank1985} (for details, see Section \ref{sec:TR-BDF2}). The resulting nonlinear systems at each time step are solved via Newton iteration and CUBLAS \cite{cuda11}, where 3-5 iterates are in general sufficient to reach the desired accuracy. The mechano-chemical system is integrated in time on the interval $[0,100]$. Due to the adaptiveness of the time-stepping method, larger time steps can be used as the system approaches the stationary state. The time stepping is terminated earlier if the $L^{\infty}$-norm of the time derivative is less than $10^{-8}$ at some point.

Fig.~\ref{Fig:6} illustrates the computational performance of the GPU-based numerical solvers, developed for this work. We measured computational times on a single Nvidia GTX 980 graphics card with 16 streaming multiprocessors (SM) and 2048 CUDA cores. The implementation of explicit Euler scheme executes simulations in blocks, each consisting of 1024 threads. A single block processes one single simulation, thus each thread processes four pixels in a spatial grid. The whole model integration process is done inside one kernel call. This makes it possible to use the shared memory of the device (a very limited amount of extremely fast block-level cache memory) to store intermediate results, but it does result in a high register usage. Thus, a single SM can run only one block instead of two. Execution time of the code does not change until all 16 SM are occupied, and then grows proportionally to the number of parallel simulations. Due to the high performance of the shared memory, this approach is very efficient for small numbers of parallel simulations. 

The parallel implementation of the ROCK2 method uses a different strategy and executes vector operations by using separate lightweight kernels, which store the intermediate data in global memory of the device and consume fewer registers. This allows us to execute more blocks per SM and to achieve better occupation of the hardware. Additionally, this implementation is not limited by the requirement to fit the state vector of the system into the shared memory of the device, which is limited by 48 kilobytes in the case of Nvidia GTX 980. However, 
the synchronization of vector operations between trajectories requires additional data transfer between host and device, which results in a significant overhead when the number $n_{\text{sim}}$ of simulations is small, but becomes much less significant when $n_{\text{sim}}$ grows. This trade-off is made possible by the dramatic reduction of right hand side evaluations achieved by ROCK2 due to its improved stability properties. As a consequence, ROCK2 scales much better with respect to the number of simulations, at least for $n_{\text{sim}} \le 256$ simulations; the execution time grows sublinearly with  $O(n_{\text{sim}}^{1/3})$.

Performance of the implicit TR-BDF2 method is governed by scaling properties of functions from CUBLAS library, which dominates the overall cost. The scaling seems to be sublinear in the number of simulations up to 4096 simulations, growing roughly as $O(n_{\text{sim}}^{1/6})$ for $n_{\text{sim}} \le 256$. The code is available 
at \url{https://github.com/AlexeyKazarnikov/CILNumericalCode}. The GPU-based code has been developed with CUDA 11.0 and compiled with Microsoft Visual C++ 2019 under Microsoft Windows and GCC 9.3.0 under Ubuntu Linux.

\subsection{Multi-feature identification for large data sets}

To illustrate the impact of the multi-feature modification on the accuracy of the parameter estimation, we employ first the original CIL algorithm \cite{kazarnikov2020} in Algorithm \ref{alg:CIL1}, but using different norms to quantify the distances. We consider each model for a fixed control parameter vector $\bm{\theta}_0$, and generate a training set of $N_\text{set}$ data patterns. For the reaction-diffusion systems \eqref{introduction:0}, \eqref{gmSystem}, and \eqref{bzSystem}, we consider the least accurate case studied in \cite{kazarnikov2020} with $N_\text{set}=50$. Due to a much richer pattern variability for fixed parameter values, we employ a larger dataset of size $N_\text{set} = 200$ for the mechanochemical model \eqref{eq:mcone} and for the reaction-diffusion-ODE model \eqref{eq:rdode}. Recall that all relevant parameters for the numerical experiments are given in Table \ref{table:exp1}. 

Using the generated patterns as input data, we start by first computing the original CIL using Algorithm \ref{alg:CIL2bst} (with bootstrapping) for each of the norms in \eqref{norm:L2}-\eqref{norm:W1infprime} separately. In all cases, we use resampling, as the amount of data is not sufficient for the basic Algorithm \ref{alg:CIL2}. Then, using all norms together in Algorithm \ref{alg:CIL2bst} we compute one multi-feature vector that contains the combined information and estimate the MCIL.  For scaled data, only the three $L_2$-based norms \eqref{norm:L2}, \eqref{norm:W12prime}, \eqref{norm:W12} are used, since min-max normalisation removes the information about absolute values required for the maximum-based norms \eqref{norm:Linf}, \eqref{norm:W1inf}, \eqref{norm:W1infprime}. In addition, for far-from-equilibrium patterns, produced by the reaction-diffusion-ODE model \eqref{eq:rdode}, gradient-based norms are computed with respect to the diffusive component $w(x,t)$ only, as the non-diffusive component $u(x,t)$ is non-differentiable.

For each model and for each of the considered settings considered, we employ an adaptive Metropolis-Hastings algorithm \cite{haario:1} to sample a chain of 12000 parameter values from the posterior distribution to determine parameter estimates and to visualise and quantify the remaining uncertainty. 
The results of this experiment are summarised in Fig.~\ref{Fig:1}-\ref{Fig:3} below. Fig.~\ref{Fig:1} shows the posterior distributions for individual distances together with the regions obtained with the MCIL approach, for all models. It can be concluded that different norms ''highlight'' different regions in parameter space, while the posterior distribution for multi-feature CIL vector (denoted by yellow colour) corresponds well to the region of intersection. Fig.~\ref{Fig:2}-\ref{Fig:3} illustrate the improvement of accuracy for non-scaled and scaled data respectively, showing the multi-feature posteriors together with the $L_2$-based posterior and patterns obtained for four verification values in parameter space. It can be seen clearly from those figures that applying a multi-feature approach significantly increases the accuracy of the identification for the same, fixed data set. 

\begin{figure}[t]
	\centering
	\includegraphics[width=0.99\textwidth]{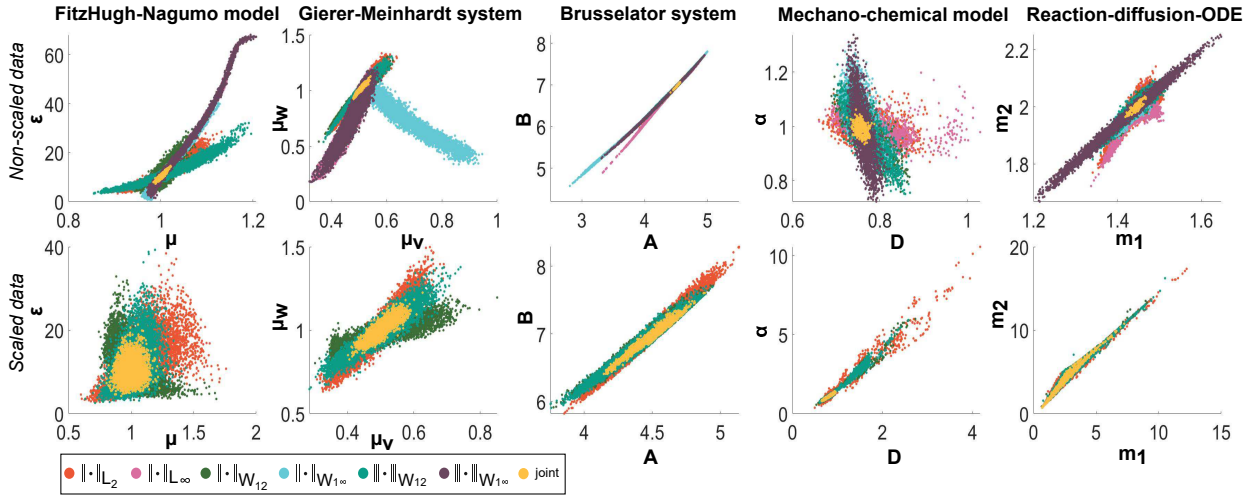}
	\caption{Posterior distributions of the parameters obtained by applying the single-feature and multi-feature CIL approach to the same training set of $N_\text{set}$ patterns (using $N_\text{set} = 50$ and $N_\text{set} = 200$ for the reaction-diffusion equations and for the other two examples, respectively).
	\label{Fig:1}} 
\end{figure}

\begin{figure}[t]
	\centering
	\includegraphics[width=0.99\textwidth]{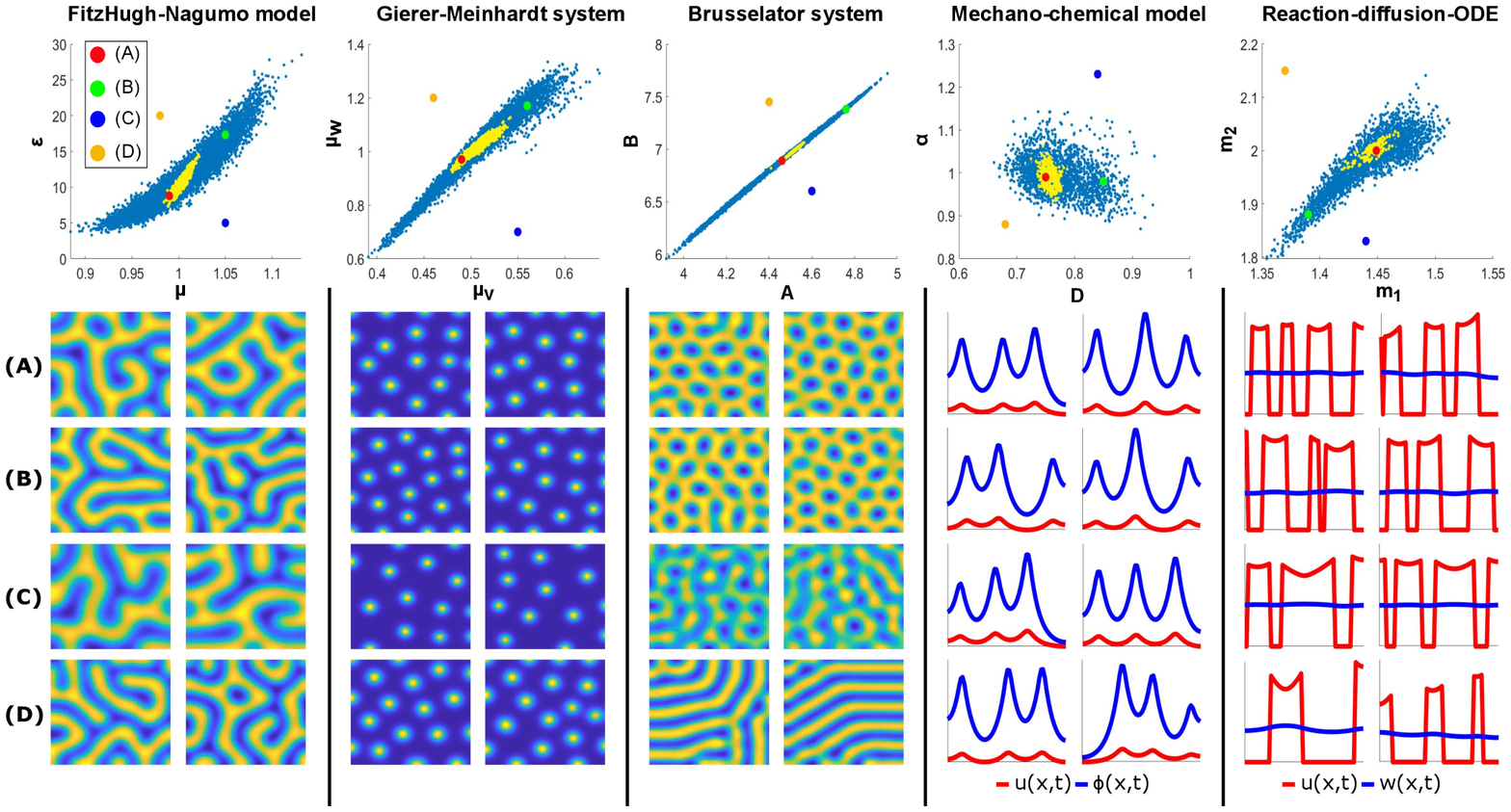}
	\caption{Top: Comparison of the posterior for the single-feature CIL approach based on the $L_2$ norm (blue) and for the multi-feature CIL approach (yellow) for non-scaled pattern data ($N_\text{set}$ as specified in Tab.~\ref{table:exp1}). Bottom: Visual inspection of the identification accuracy at a number of verification values of the control parameters, chosen within and outside the estimated posterior regions. 
	}
	\label{Fig:2}
\end{figure}

\begin{figure}[t]
	\centering
	\includegraphics[width=0.99\textwidth]{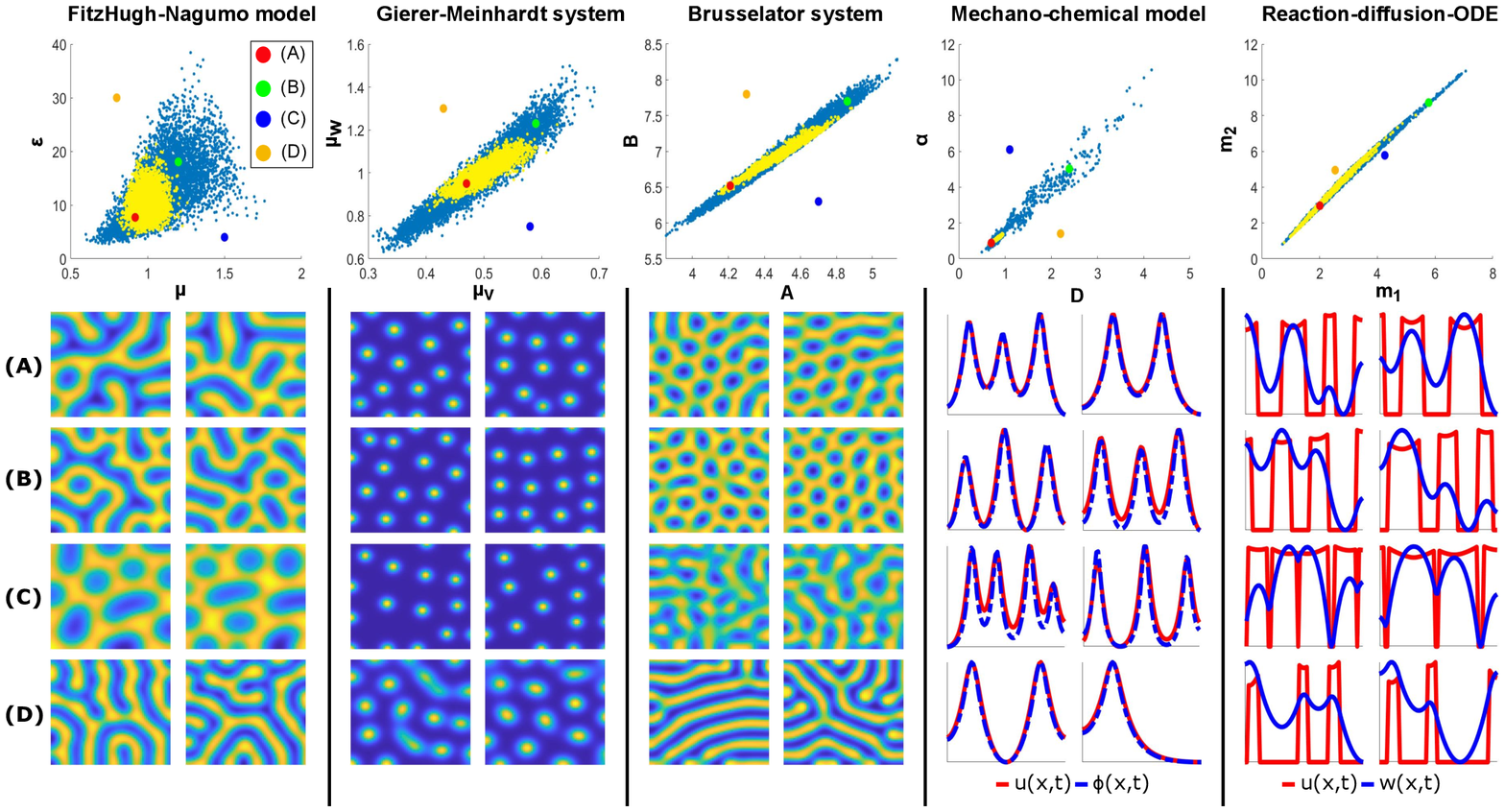}
	\caption{Top: Comparison of the posterior for the single-feature CIL approach based on the $L_2$ norm (blue) and for the multi-feature CIL approach (yellow) for scaled pattern data ($N_\text{set}$ as specified in Tab.~\ref{table:exp1}). Bottom: Visual inspection of the identification accuracy at a number of verification values of the control parameters, chosen within and outside the estimated posterior regions.
	}
	\label{Fig:3}
\end{figure}

Quantification of distances between patterns depends on a chosen norm. While mapping from high-dimensional pattern data to scalar numbers results in a loss of information, different norms emphasise different features of the data. Therefore, taking into account a set of norms we obtain the 'intersection' of the respective posterior regions, see Fig.~\ref{Fig:1}--\ref{Fig:3}. For some models and some norms, the posteriors coincide or are very similar, but for others the difference can be significant. As stated above, the overhead of computing multiple norms and of dealing with larger generalised correlation integral vectors is inexpensive compared to carrying out the model simulations. Thus, it always pays off to consider all admissible norms.

 
For the reaction-diffusion models and a data set of size $N_\text{set}=50$, multi-feature identification achieves a much more accurate detection of changes in model parameters than the performance previously achieved with the $L_2$-based CIL approach in \cite{kazarnikov2020}, which had been comparable to the performance of the 'naked eye'. Similar observations hold for the mechano-chemical model with one diffusing morphogen and a data set of size $N_\text{set}=200$ in the case of non-scaled data. The CIL approach with one norm works reasonably well, and the MCIL modification significantly improves the accuracy. However, considering min-max scaled data, the resulting posteriors with a single norm are significantly wider. This may be explained by the fact that the tissue and morphogen patterns have identical profiles. Thus, scaling removes the information about absolute values of each component and their ratio (see the respective pictures in Fig.~\ref{Fig:2}-\ref{Fig:3}), forcing the algorithm to work with highly changeable one-dimensional pattern curves. In this case, the derivative-based Sobolev norms \eqref{norm:W12} and \eqref{norm:W12prime} are very helpful in the detection of small changes in curve variability, thus significantly reducing the width of the posterior distribution for the MCIL, see Fig.~\ref{Fig:1}.

\subsection{Synthetic likelihood for limited data}

We examine the performance of the synthetic likelihood idea. First, we consider the same training set of $N_{set}$ patterns as used above, and define the cost function $f(\bm{\theta})$, by using Algorithm~\ref{alg:CIL3bst} with norms \eqref{norm:L2}-\eqref{norm:W1infprime}. Here we again use resampling to avoid large computational times. For each evaluation of the cost function $f(\bm{\theta})$, we simulate a set of $N_\text{syn} = 1000$ patterns and compute distances between them by using different norms. 
Next, by using bootstrapping, we sub-sample a set of $n_\text{CIL}=1000$ correlation integral vectors at $\bm{\theta}$ to obtain reliable estimates for the mean vector $\bm{\mu}_{\bm{\theta}}$ and for the covariance matrix $\bm{\Sigma}_{\bm{\theta}}$ of the multi-feature SCIL. Finally, we compute the output of the cost function as shown in Algorithm \ref{alg:CIL3bst}. 

\begin{figure}[t]
	\centering
	\includegraphics[width=0.99\textwidth]{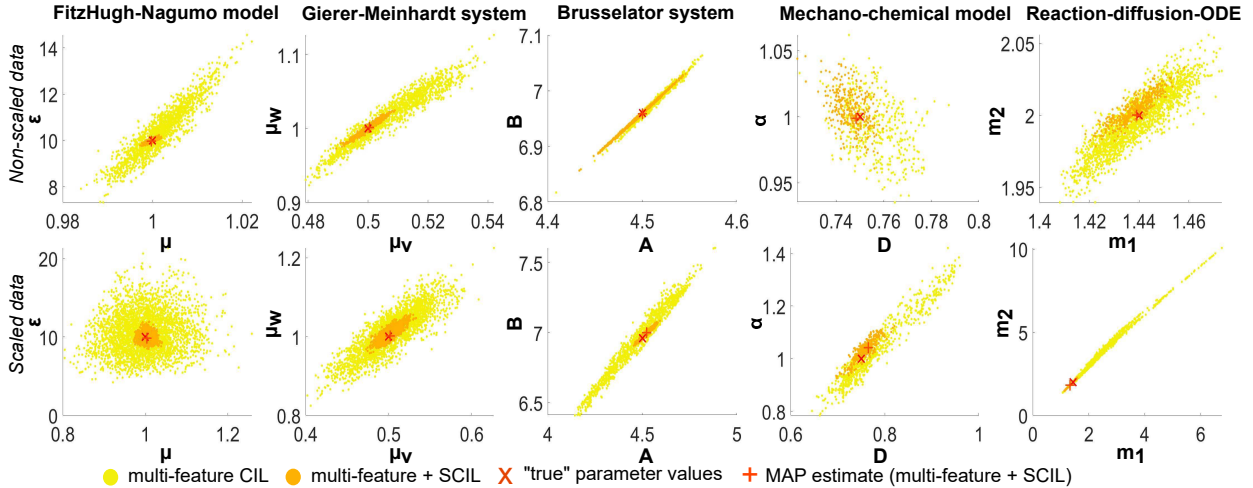}
	\caption{Posterior distributions of the parameters obtained by applying the multi-feature CIL and the synthetic (multi-feature) CIL approach to the same training set of $N_\text{set}$ patterns (with $N_\text{set}$ as specified in Tab.~\ref{table:exp1}).}
	\label{Fig:4}
\end{figure}

\begin{figure}[t]
	\centering
	\includegraphics[width=0.99\textwidth]{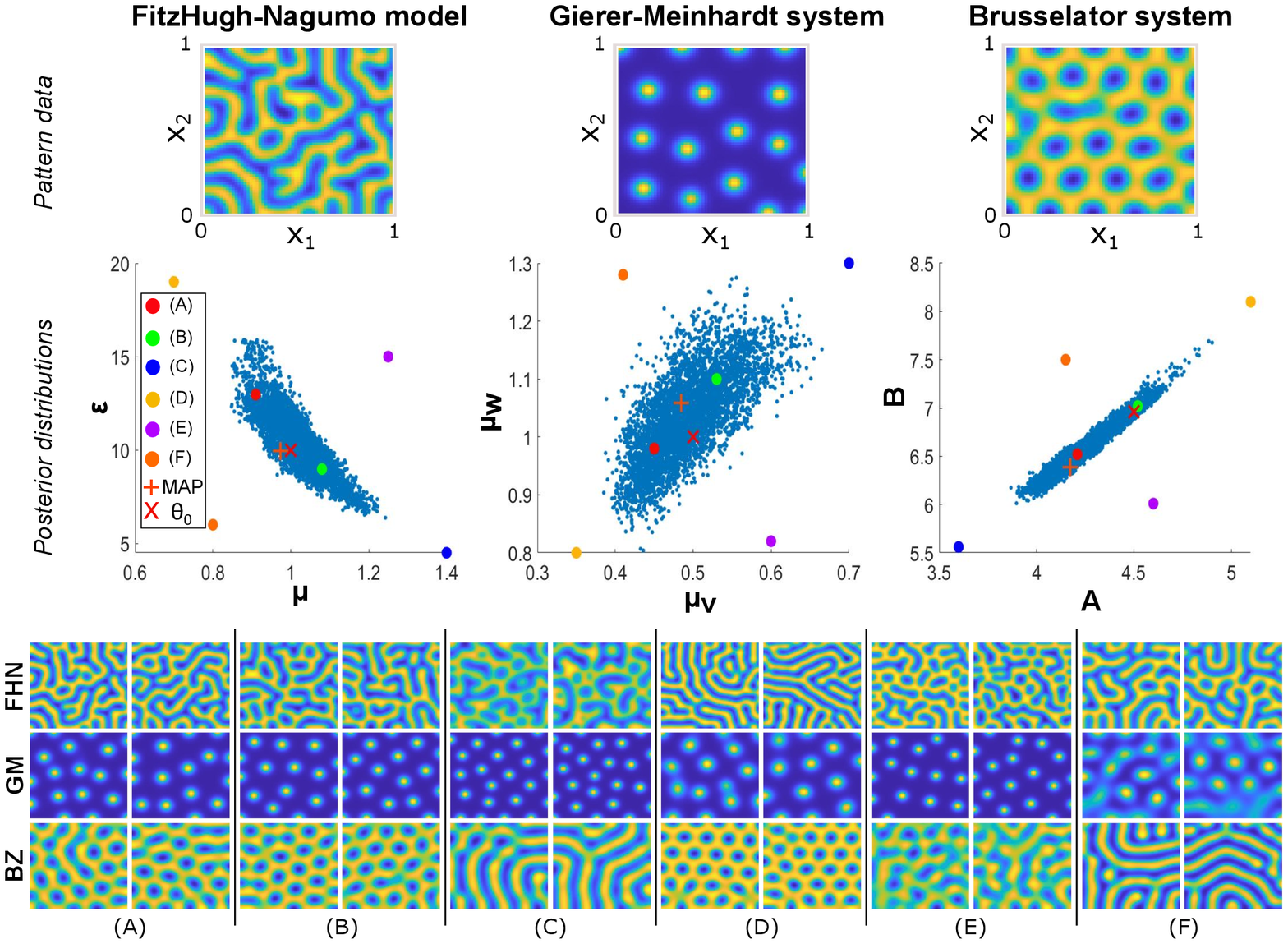}
	\caption{Posterior distributions of the parameters for the reaction-diffusion models, obtained by applying the synthetic (multi-feature) CIL with a single min-max scaled pattern. }
	\label{Fig:5}
\end{figure}
To study the performance of this cost function, we use an adaptive Metropolis-Hastings algorithm to sample a chain of 12000 parameter values for each of the models. The results are shown in Fig.~\ref{Fig:4}. The synthetic likelihood-based cost function 
clearly seems to further increase the accuracy of the multi-feature identification results without any resizing of the data set. 
Moreover, due to the fact that the SCIL is generated now on the fly by model-generated data independently of the size of the  training data set, the applicability of the approach and the quality of the identification process is much less affected by the amount of patterns $N_\text{set}$. For each model under study, we also examined the behaviour of the method for smaller training sets of size, e.g.,  $N_\text{set}/2$ and $N_\text{set}/4$. The performance of the method was satisfactory in all these cases: posterior distributions were naturally getting larger when decreasing the number of patterns, but at the same time they remained clearly bounded.


Finally, let us discuss the situation for the most limited case of information, consisting of only one min-max scaled pattern with about 5 'wavelengths' as data. Here, we restrict ourselves to the reaction-diffusion models. Naturally, the amount of information encoded in the picture plays a crucial role. The spatial extent of the pattern must be reasonably large to encapsulate the typical pattern variability. For the Gierer-Meinhardt and Brusselator models we use the same settings as earlier (see Fig.~\ref{Fig:5}, top), which are also comparable to those in the article \cite{schnorr2021}, while each individual pattern in~\cite{zhu2022} contains significantly more 'wavelengths' and thus more information per pattern. The patterns in the FitzHugh-Nagumo model are more rich, and we must ensure that the single data pattern contains all typical features.
We achieve this by slightly changing the diffusion coefficients, multiplying both by $0.35$, in order to change the relative spatial scale of the patterns. Otherwise, we use exactly the same settings as before.
 The results are given in Fig.~\ref{Fig:5}. Again, the posterior distributions appear to remain bounded; the variance is reasonable and the patterns at the verification values of the control parameters seem to correctly represent the 
 features of the pattern data when overlapped by the posterior, while visually differing outside the empirical distribution. The performance appears to be roughly comparable with that of the naked eye in this case.



\section{Discussion}
\label{sec:discussion}

In this paper, we have considered the problem of parameter identification by pattern data. We have introduced modifications to the Correlation Integral Likelihood approach, which allow significantly increasing the accuracy of parameter identification with the same amount of training data and extend the applicability of the method to the situation when the size of the training dataset is severely limited. The approach was tested with three classes of pattern formation models: reaction-diffusion systems exhibiting the formation of Turing patterns, mechano-chemical models producing stationary biomechanical patterns and reaction-diffusion-ODE systems leading to {\it far-from equilibrium} patterns with jump-discontinuity. The performance of the method was highly satisfactory for all considered cases. In addition, we have included experiments with min-max normalised data, which corresponds to the situation when only greyscale images of the patterns are available. For all models under study, we have provided GPU-based parallel implementations of the respective numerical algorithms.

As shown in \cite{kazarnikov2020}, the CIL method does not depend on the accuracy of the initial choice of the unknown parameter values. The algorithm also converges when starting with initial parameter values that are distant from the real ones. Thus, the CIL approach provides a robust technique for parameter identification by pattern data, even for severely limited data sets. In the numerical experiments, conducted in this work, we used known values of control parameters as starting points for MCMC sampling. However, while examining the performance of the proposed synthetic likelihood idea, we additionally verified the convergence of the method by applying the differential evolution algorithm with the same initial choices of the parameters as considered in \cite{kazarnikov2020}.

Obtaining point estimates of model parameters with the CIL approach requires a moderate number of model evaluations. In our test cases with two control parameters, the DE algorithm usually converges to the local minimum after $50$ iterations with $20$ candidates in a population.  This results in $1000 \times N_{set}$ and $1000 \times N_{syn}$ model evaluations for the MCIL and SCIL cases, respectively. In creating the parameter posteriors by MCMC methods, the amount of required model evaluations may be higher, around $N_{chain} \times N_{set}$ and $N_{chain} \times N_{syn}$, where generally $N_{chain} \geq 12000$. However, random-walk-based MCMC are conceivable 
by using suitable techniques, such as multi-level MCMC methods \cite{dodwell2015} or local approximation methods for the posterior \cite{springer2021,conrad2018}. It results in a further improvement of the method. In our examples, we have employed only parameter vectors of dimension two, but the approach scales with respect to the parameter dimension in the same way as optimisation and the MCMC sampling typically do. 

The accuracy of parameter identification achieved by the CIL approach can be improved by introducing additional mappings from pattern data to scalar values in the CIL scheme. In this work, we defined these mappings as different Lebesgue and Sobolev norms, but additional norms, such as total variation norm or bounded Lipschitz distance may be helpful, depending on the model. Additionally, different pattern representations, such as resistance distance histograms, introduced in \cite{schnorr2021},  can be used as features. We leave this topic for further investigation.

The alternative approach to identify model parameters from pattern data using artificial neural networks, discussed in \cite{zhu2022,schnorr2021}, may be more expensive to train, but the evaluation of the trained model is very fast. Moreover, once trained, the model may be applied to different sets of pattern data, provided that the unknown parameter values (or nearby values) were shown to the neural network during training. The accuracy of the approach can be further improved by increasing the amount of training data or by adjusting the architecture of the network. However, for parameter values that are out-of-distribution and have not been used during training, the performance can be very poor and unpredictable. Typically,
a re-training of the model in the new parameter region is required, which results in additional computational cost. Moreover, the required amount of training data rapidly grows with the number of control parameters \cite{zhu2022}, which might become prohibitively large for practical applications. Finally, reliable uncertainty quantification is not available. 

The two approaches could be combined to increase the efficiency of parameter identification. First, the CIL approach can be used to identify the parameter values that are sufficiently accurate and localise the region of interest in the parameter space by constructing a posterior distribution. The information on the sets of parameters and the respective patterns, which was gained during the estimation, can be further used to create an optimal training set for a suitable neural network model. Once constructed, the trained network can be used to obtain fast estimates of model parameters by pattern data. We leave this topic for further studies.

{ Another important aspect of model selection is related to the choice of initial conditions. In this work, the focus lay on models for {\it de novo} pattern formation and hence the initial data were assumed to be small perturbations of spatially homogeneous steady states. However, in many applications it is necessary to deal with tightly constrained initial 
conditions in order to obtain specific types of patterns, with a prescribed number of 'elements' or 'features' in specific locations. It can be shown numerically that in all considered pattern formation models, strong, localised perturbations may lead to pre-selection of a particular family of patterns, e.g. formation of patterns organised around a pre-selected spot. 
Identifying such patterns, which are not {\it de novo}, requires including the choice of the initial condition in the modelling process, e.g., by accounting for the experimental conditions generating the observed patterns, which is also possible in our framework. Nevertheless, as solving this problem is an application-dependent topic, we leave it for future work.}

\subsection*{Acknowledgements}

This work is supported by the Deutsche Forschungsgemeinschaft (German Research Foundation) under Germany’s Excellence Strategy
EXC 2181/1 - 390900948 (the Heidelberg STRUCTURES Excellence Cluster), by the Collaborative Research Center (SFB) 1324, and by the Centre of Excellence of Inverse Modelling and Imaging (CoE), Academy of Finland,  decision number 312 125. { Additionally, we would like to thank Andrew L. Krause for his valuable comments and stimulating suggestions.}

{\small
\bibliographystyle{ieeetr}
\bibliography{paper_correction_text.bib}
}

\pagebreak

\begin{appendix}

\begin{center}
    \Large
    \textbf{Supplementary Materials}
\end{center}

\section{Additional algorithms: MCIL and SCIL with bootstrapping}
\label{sec:algos}

In Algorithm \ref{alg:CIL2bst}, a bootstrapping resampling procedure to improve the parameter estimates for the multi-feature CIL approach in the limited data case is described. Similarly, Algorithm \ref{alg:CIL3bst} provides details of bootstrapping for the synthetic CIL approach.

\begin{algorithm}[h!]
	\SetKwInOut{Input}{Input}\SetKwInOut{Output}{Output}
	\SetKwInOut{Data}{Data}
	\SetKwInOut{Result}{Result}
	\Data{A set $\bm{s}_\text{data}$ of discretised patterns $\bm{s}_1,\ldots,\bm{s}_{N_\text{set}}$, with unknown model parameter $\bm{\theta}_0$}
	\Input{Family of norms $\|\cdot\|_\alpha$, $\alpha=1,\ldots,N_\text{dist}$}
	\Input{$n_\text{CIL} \ \ \ldots$ \ \ Number of generalised correlation integral vectors to estimate MCIL}
	\Result{Parameters $\bm{\mu}_0$ and $\bm{\Sigma}_0$ for the cost function $f(\bm{\theta})$ in \eqref{cil:CF} in the case of the MCIL}
	\Begin
	{
		1. Divide $\bm{s}_\text{data}$ into two subsets $\widetilde{\bm{s}}^1$ and $\widetilde{\bm{s}}^2$, each containing $N_\text{set}/2$ patterns
		
		2. \For{$k=1,\ldots,n_\mathrm{CIL}$}
		{
		2.1 Construct two sets $\bm{s}^1$ and $\bm{s}^2$ by randomly selecting with replacement $N_\text{set}/2$ patterns from the sets $\bm{\tilde{s}}^1$ and $\bm{\tilde{s}}^2$ respectively
		
		2.2 Initialise the $k$th correlation integral vector $\bm{y}^{k}$ to be the empty vector

		2.3 \For{$\alpha=1,\ldots,N_\mathrm{dist}$}
		{
			* Using $\|\cdot\|_\alpha$ compute the distances $\|\bm{s}^1_{i}-\bm{s}^2_{j}\|_\alpha$ between all patterns in $\bm{s}^1$ and $\bm{s}^2$
			
			* Using $\|\cdot\|_\alpha$ compute the part $\bm{y}^{k}_\alpha$ of the gen.~correlation integral vector via \eqref{eq:CIL}
			
			* Concatenate the current vector $\bm{y}^{k}$ and $\bm{y}^{k}_\alpha$
		}
        }
	
    3. Using the samples $\bm{y}^{k}$ computed in Step 2, estimate mean $\bm{\mu}_0$ and covariance matrix $\bm{\Sigma}_0$ of the (multi-feature) generalised correlation integral vector $\bm{y}_0$ of the training data

	}
	\caption{Using bootstrapping to compute the cost function $f(\bm{\theta})$ in the MCIL approach\!\!\!}
	\label{alg:CIL2bst}
\end{algorithm}

\begin{algorithm}[t]
	\SetKwInOut{Input}{Input}\SetKwInOut{Output}{Output}
	\SetKwInOut{Data}{Data}
	\SetKwInOut{Result}{Result}
	\Data{A set $\bm{s}_\text{data}$ of discretised patterns $\bm{s}_1,\ldots,\bm{s}_{N_\text{set}}$, with unknown model parameter $\bm{\theta}_0$}
	\Input{$\|\cdot\|_\alpha  \ \ldots$ \ a family of norms with $\alpha=1,\ldots,N_\text{dist}$}
	\Input{$\bm{\theta} \quad\ \ \ \ldots$ \ the parameter value where the cost function should be evaluated}
	\Input{$N_\text{syn} \ \ \ldots$ \ the number of model-generated (synthetic) patterns at $\bm{\theta}$}
	\Input{$n_\text{CIL} \ \ \ldots$ \ \ Number of generalised correlation integral vectors to estimate MCIL}
	\Output{$f(\bm{\theta}) \ \ \ldots$ \ the value of the  cost function at $\bm{\theta}$}	
	\Begin
	{
		1. Simulate $N_\text{syn}$ patterns $\bm{s}_\text{syn} := \{\bm{s}_i(\bm{\theta}):i=1,\ldots,N_\text{syn}\}$, with parameter $\bm{\theta}$
		
		2. \For{$k=1,\ldots,n_\mathrm{CIL}$}
		{
		2.1 Construct $\bm{s}^1$ by randomly selecting with replacement $N_\text{set}$ patterns from $\bm{s}_\text{syn}$
		
		2.2 Construct $\bm{s}^2$ by randomly selecting with replacement $\widetilde{N} = N_\text{syn} - N_\text{set}$ patterns from the remaining patterns in $\bm{s}_\text{syn}$
		
		2.3 Initialise the correlation integral vector $\bm{y}_{\bm{\theta}}^{k}$ to be the empty vector

		2.4 \For{$\alpha=1,\ldots,N_\mathrm{dist}$}
		    {
			* Compute the distances $\|\bm{s}^{1}_{i}-\bm{s}^{2}_{j}\|_\alpha$ between all patterns in $\bm{s}^{1}$ and $\bm{s}^{2}$ with $\|\cdot\|_\alpha$
			
			* Compute part of correlation integral vector $\bm{y}^{k}_{\bm{\theta},\alpha}$ as defined in \eqref{eq:CIL_synth_kl} using $\|\cdot\|_\alpha$
			
			* Concatenate the current vector $\bm{y}_{\bm{\theta}}^{k}$ and $\bm{y}^{k}_{\bm{\theta},\alpha}$
		    }
        }		
        3. Using the samples $\bm{y}_{\bm{\theta}}^{k}$ computed in Step 2, estimate mean $\bm{\mu}_{\bm{\theta}}$ and covariance matrix $\bm{\Sigma}_{\bm{\theta}}$ of the (multi-feature) correlation integral vector $\bm{y}_{\bm{\theta}}$ of the synthetic data at $\bm{\theta}$

		4. Randomly select a subset $\bm{s}_{\bm{\theta}}^{2}$ of $\widetilde{N}$ patterns from $\bm{s}_\text{syn}$

		5. Using the computed estimates $\bm{\mu}_{\bm{\theta}}$ and $\bm{\Sigma}_{\bm{\theta}}$ and $\widetilde{\bm{y}}(\theta)$ as defined in~\eqref{eq:CIL_synth} (with $\bm{s}_{\bm{\theta}}^{2}$ instead of $\bm{s}_{\bm{\theta}}^{k,2}$), compute the value $f(\bm{\theta})$ of the cost function at $\bm{\theta}$ from \eqref{cil:CF_synth}
	}
	\caption{Using bootstrapping to evaluate $f(\bm{\theta})$ for SCIL at a single parameter value $\bm{\theta}$\!\!\!}
	\label{alg:CIL3bst}
\end{algorithm}

\section{Numerical discretisation of the studied models}
\label{sec:discretise}


\subsection{Spatial discretisation using the Method of Lines}
To discretise the PDE models in space, we apply the Method of Lines (MoL) and finite difference approximations. Let us discuss the discretisation process for each system separately.

\subsubsection{Reaction-diffusion systems}
 The square domain $\Omega$ is discretised using a uniform grid with fixed step size $h=1/(M_\text{dim}-1)$, $M_\text{dim}\in \mathbb{N}$ leading to a finite set of grid points: 
\[
\big\{\bm{x}_{i,j}=((i-1)h,(j-1)h),\,i,j=1,\dots,M_\text{dim}\big\}.
\]

To reduce the reaction-diffusion system \eqref{eq:rdgen} to a finite set of $2M_\text{dim}^2$ ordinary differential equations (ODEs), we define for each grid point the time-dependent functions
\[
v_{i,j}(t) = v(\bm{x}_{i,j},t),
\ \ 
w_{i,j}(t) = w(\bm{x}_{i,j},t),
\quad
i,j=1,\ldots,M_\text{dim},
\]
and discretise the Laplace operator by the five-point stencil \cite{grossman:1,hupkes:1}
\[
\Delta v\approx\frac{v_{i+1,j}+v_{i-1,j}+v_{i,j+1}+v_{i,j-1}-4v_{i,j}}{h^{2}}=\nabla_{h}^{2}v_{i,j}.
\]
Substituting the approximations leads to the system:
\begin{eqnarray}
	{\textstyle \frac{\text{d}}{\text{d}t}} v_{i,j}(t) & = & \nu_1\nabla_{h}^{2}v_{i,j}(t)+f(v_{i,j}(t),w_{i,j}(t)), \nonumber \\[1ex] 
	{\textstyle \frac{\text{d}}{\text{d}t}} w_{i,j}(t) & = & \nu_2\nabla_{h}^{2}w_{i,j}(t)+g(v_{i,j}(t),w_{i,j}(t)), \quad \label{MOLSystem}
	i,j=1,2,\dots,M_{dim}.
\end{eqnarray}
The Neumann boundary conditions are taken into account by using a one-sided first-order difference scheme \cite{grossman:1}, which leads to the following conditions for the 'ghost' values $\bm{x}_{0,j}$, $\bm{x}_{N+1,j}$, $\bm{x}_{i,0}$, $\bm{x}_{i,N+1}$:
\begin{eqnarray*}
	& v_{0,j}(t) \equiv v_{1,j}(t),
	\,
	v_{M_{dim}+1,j}(t) \equiv v_{M_{dim},j}(t), \\
	& w_{0,j}(t) \equiv w_{1,j}(t),
	\,
	w_{M_{dim}+1,j}(t) \equiv w_{M_{dim},j}(t), \\
	& v_{i,0}(t) \equiv v_{i,1}(t),
	\,
	v_{i,M_{dim}+1}(t) \equiv v_{i,M_{dim}}(t), \\
	& w_{i,0}(t) \equiv w_{i,1}(t),
	\,
	w_{i,M_{dim}+1}(t) \equiv w_{i,M_{dim}}(t).
\end{eqnarray*}

\subsubsection{Mechano-chemical models}
Before discretising the non-linear mechanochemical models with the MoL, it is necessary to obtain analytical expressions for spatial derivatives in the respective equations. In \eqref{eq:morphogen}, the diffusion operator for the morphogen $\phi(x,t)$ can be written in the form:
\[
\dfrac{\text{d}}{\text{d}x}\left(\dfrac{\phi'(x,t)}{\sqrt{1+u'(x,t)^2}}\right) = \dfrac{\phi''(x,t)}{\sqrt{1 + u'(x,t)^2}} - \dfrac{\phi'(x,t) u'(x,t) u''(x,t)}{\left( 1 + u'(x,t)^2 \right)^{\frac{3}{2}}}.
\]
The respective term in the equation for tissue movement \eqref{eq:tissue} can be expanded as:
\begin{equation}
\dfrac{\text{d}}{\text{d}x}\left(\dfrac{\kappa'(u)-\bar{\kappa}'(\phi)}{\sqrt{1+u'(x,t)^2}}\right)
=
\dfrac{\kappa''(u)-\bar{\kappa}''(\phi)}{\sqrt{1 + u'(x,t)^2}}
-
\dfrac{\left( \kappa'(u) - \bar{\kappa}'(\phi) \right) u'(x,t) u''(x,t)}{\left( 1 + u'(x,t)^2 \right)^{\frac{3}{2}}}.
\label{eq:tissue_curvature}
\end{equation}
Equation \eqref{eq:tissue_curvature} contains the spatial derivatives $\kappa'(u)$ and $\kappa''(u)$ of the tissue curvature \eqref{eq:curvature}. These terms can be written as:
\begin{align}
\kappa'(u)
&=
\dfrac{u'''(x,t)}{(1+u'(x,t)^2)^{\frac{3}{2}}}
-
3\,\dfrac{ u'(x,t) u''(x,t)^2}{(1+u'(x,t)^2)^{\frac{5}{2}}},
\label{eq:curvature_deriv1}	\\[0.5ex]
\kappa''(u) 
&=
\dfrac{u^{(iv)}(x,t)}{(1+u'(x,t)^2)^{\frac{3}{2}}}
-
3\,\dfrac{ u'(x,t) u''(x,t) u'''(x,t)}{(1+u'(x,t)^2)^{\frac{5}{2}}} 
- \vspace{2mm} \nonumber\\[0.5ex]
&- 3\,\dfrac{   u''(x,t)^3 + 2 u'(x,t) u''(x,t) u'''(x,t) }{(1+u'(x,t)^2)^{\frac{5}{2}}}
+
15\,\dfrac{  u'(x,t)^2 u''(x,t)^3 }{(1+u'(x,t)^2)^{\frac{7}{2}}}.
\label{eq:curvature_deriv2}	
\end{align}

Next, the one-dimensional interval $\Omega = (0,L)$ is discretised using a uniform grid with fixed step size $h = L/(M_\text{dim}-1)$, $M_\text{dim} \in \mathbb{N}$, leading again to a finite set of grid points:
\[
\lbrace x_j = (j-1)h:j=1,\ldots,M_\text{dim} \rbrace.
\]
We define time-dependent functions
\[
u_j(t) = u(x_j,t),
\ \ 
\phi_j(t) = \phi(x_j,t),
\quad
j=1,\ldots,M_\text{dim}.
\]
and approximate the spatial derivatives again by finite differences, in particular using the forward difference scheme to approximate the first derivatives of a function $f(x)$, i.e.,
\[
f'(x_j) \approx D_h[f](x_j) = \dfrac{f(x_{j+1})-f(x_j)}{h},
\]
and central difference schemes for higher order derivatives, i.e.,
\begin{align*}
f''(x_j) & \approx D_h^2[f](x_j) = \dfrac{f(x_{j+1})-2f(x_j)+f(x_{j-1})}{h^2},\\[0.5ex]
f'''(x_j) & \approx D_h^3[f](x_j) = \dfrac{f(x_{j+2})-2f(x_{j+1})+2f(x_{j-1})-f(x_{j-2})}{2h^3},\\[0.5ex]
f^{(iv)}(x_j) & \approx D_h^4[f](x_j) = \dfrac{f(x_{j+2})-4f(x_{j+1})+6f(x_j)-4f(x_{j-1})+f(x_{j-2})}{h^4}.
\end{align*}
This allows to derive finite difference approximations for the curvature terms in \eqref{eq:curvature}, \eqref{eq:curvature_deriv1}, \eqref{eq:curvature_deriv2}, which we will not explicitly write down:
\begin{eqnarray*}
\kappa(u(x_j,t)) & \approx : & \kappa_h(u_{j-1}(t),u_j(t),u_{j+1}(t)), \\
\kappa'(u(x_j,t)) & \approx : & \kappa^{(i)}_h(u_{j-2}(t),u_{j-1}(t),u_j(t),u_{j+1}(t),u_{j+2}(t)), \\
\kappa''(u(x_j,t))& \approx : & \kappa^{(ii)}_h(u_{j-2}(t),u_{j-1}(t),u_j(t),u_{j+1}(t),u_{j+2}(t)),
\end{eqnarray*}
as well as for the derivatives of the locally preferred curvature, i.e.,
\[
\bar{\kappa}_h^{(i)}(\phi_j(t),\phi_{j+1}(t)) = -\beta D_h[\phi](x_j),
\quad
\bar{\kappa}_h^{(ii)}(\phi_{j-1}(t),\phi_j(t),\phi_{j+1}(t)) = -\beta D^2_h[\phi](x_j).
\]

In these expressions, Neumann boundary conditions are taken into account, giving the following expressions for the 'ghost' nodes:
\[
\begin{array}{c}
u_0(t) = u_1(t); \, u_{N+1}(t) = u_N(t), \\
u_{-1}(t) = u_3(t) - 2u_2(t) + 2u_1(t); \,
u_{N+2}(t) = 2u_N(t) - 2u_{N-1}(t) + u_{N-2}(t), \\
{\phi}_0(t) = {\phi}_1(t); \,
{\phi}_{N+1}(t) = {\phi}_N(t).
\end{array}
\]

Thus, one arrives at the following finite difference approximation for the mechano-chemical model with one diffusing morphogen, 
for $j=1,\ldots,M_\text{dim}$:
\begin{equation}
	\begin{array}{l}
		\dfrac{\text{d}\phi_j(t)}{\text{d}t} 
		=
		D\left( 
		\dfrac{D^2_h[\phi](x_j)}{\sqrt{1 + D_h[u](x_j)^2}} 
		- 
		\dfrac{D_h[\phi](x_j) D_h[u](x_j) D^2_h[u](x_j)}{\left( 1 + D_h[u](x_j)^2 \right)^{\frac{3}{2}}}
		\right)
		-
		\alpha\phi_j(t) 
		+ 
		f(\kappa_h),
		\vspace{2mm}\\
		\tau\dfrac{\text{d}u_j(t)}{\text{d}t} = -L
		\Bigg[
		\dfrac{\kappa^{(ii)}_h-\bar{\kappa}^{(ii)}_h}{1 + D_h[u](x_j)^2}
		-
		\dfrac{\Big( \kappa^{(i)}_h - \bar{\kappa}^{(i)}_h \Big) D_h[u](x_j) D^2_h[u](x_j)}{\left( 1 + D_h[u](x_j)^2 \right)^2}
		-
		\left(\kappa_h-\bar{\kappa}_h\right)\kappa_h^2 + \lambda\kappa_h
		\Bigg].
	\end{array}
	\label{eq:mc:MOL}
\end{equation}
Finally, the global arc length integral can be approximated by the trapezoidal rule to give
\begin{equation}
S(u) 
\approx  
h
\left(
	\sqrt{1+D_h[u](x_1)^2)}
		+
	2\sum\limits_{j=2}^{M_\text{dim}-1}\sqrt{1+D_h[u](x_j)^2)}
		+
	\sqrt{1+D_h[u](x_{M_\text{dim}})^2}
\right). 
\label{eq:mc:arc}
\end{equation}


\subsection{Time integration of the equations under study}

To finally map the problems onto a GPU, we also need to discretise the ODE models in Section 9.1 with respect to time. This procedure is discussed below.

\subsubsection{Reaction-diffusion systems}
\label{sec:ROCK2}


To simulate the discretised reaction-diffusion systems \eqref{MOLSystem} we use the ROCK2 numerical method, which was introduced in \cite{abdulle2001}. It is an explicit, second-order stabilised Runge-Kutta method that possesses a large stability domain along the negative real axis, obtained by approximating the optimal stability polynomials of second order by suitable orthogonal  polynomials, and by exploiting the inherent three-term recurrence relations. It is thus especially efficient for finite-difference approximations of parabolic PDEs obtained by the MoL \cite{hairer1996}.


Let us denote the state vector of system \eqref{MOLSystem} by $\bm{w} \in \mathbb{R}^{m}$ and the right-hand side of this system by $\bm{F}(\bm{w}) : \mathbb{R}^{m} \rightarrow \mathbb{R}^{m}$, where $m = 2M_{dim}^2$, such that
\[
{\textstyle \frac{\text{d}}{\text{d}t}} \bm{w}(t) = \bm{F}\big(\bm{w}(t)\big), \quad t > 0, \quad \text{and} \quad \bm{w}(0)=\bm{w}_0.
\] 
Denote the approximation of $\bm{w}(t_n)$ at time step $t_n$, $n \geq 1$, by $\bm{w}_n$. Then, the approximation $\bm{w}_{n+1}$ at the next time step is obtained as follows:
\begin{align}
		\bm{g}_0 &= \bm{w}_n, \nonumber\\[0.5ex]
		\bm{g}_1 &= \bm{w}_n + h_n\mu_1^s \bm{F}(\bm{g}_0), \nonumber\\[0.5ex]
		\bm{g}_j &= h_n\mu_j^s \bm{F}(\bm{g}_{j-1}) - \nu_j^s \bm{g}_{j-1} - \kappa_j^s \bm{g}_{j-2}, \quad j=2,\ldots,s-2, \label{eq:3term}\\[0.5ex]
		\bm{g}_{s-1}&=\bm{g}_{s-2}+h_n\sigma^s \bm{F}(\bm{g}_{s-2}), \nonumber\\[0.5ex]
		\bm{g}_s^* &= \bm{g}_{s-1} + h_n\sigma^s \bm{F}(\bm{g}_{s-1}), \nonumber\\[0.5ex]
		\bm{w}_{n+1} &= \bm{g}_s^* - h_n\sigma^s(1-\frac{\sigma^s}{\tau^s})(\bm{F}(\bm{g}_{s-1})-\bm{F}(\bm{g}_{s-2})).\nonumber
\end{align}
Here, $h_n$ is the current time step size. The stages $j=0,1,\ldots,s-2$ are computed using the three-term recurrence relation
in \eqref{eq:3term}, while the final two stages $j=s-1,s$ are determined by a suitable two-stage finishing procedure. The number of stages $s$ is determined from the relation:
\[
h_n \rho\big(\text{\bf D} \bm{F} 
(\bm{w}_n)\big) \; \leq \; (0.9 s)^2,
\]
where $\rho(\cdot)$ denotes the spectral radius of a matrix and $\text{\bf D} \bm{F}(\bm{w}_n)$ is the Jacobian of $\bm{F}$ at $\bm{w}_n$. In numerical simulations, we approximate this spectral radius by employing Gershgorin's Circle Theorem~\cite{gershgorin1931}. For more details, see \cite{abdulle2001}. 

The error after each time step is computed by
\[
err_{n+1} = \left\lbrace \dfrac{1}{m} \sum\limits_{k=1}^m \left( \dfrac{|\bm{w}_{n+1,k} - \bm{g}_{s,k}^*|^2}{a_{tol} + r_{tol}\,\max(|\bm{w}_{n,k}|,|\bm{w}_{n+1,k}|)} \right) \right\rbrace^{\frac{1}{2}},
\]
where $a_{tol}$ and $r_{tol}$ are respectively absolute and relative tolerances, and the new time step is defined by
\[
h_{new} = \eta\, h_n \left(\dfrac{1}{err_{n+1}}\right)^{\frac{1}{2}} \dfrac{h_n}{h_{n-1}}\left(\dfrac{err_n}{err_{n+1}}\right)^{\frac{1}{2}},
\]
where we set in numerical simulations $\eta = 0.8$.

Due to the large stability domain and the fact that the method is second-order, the ROCK2 numerical method leads to a dramatically reduced number of evaluations of the right-hand side $\bm{F}$ than the explicit Euler method employed in our previous implementation \cite{kazarnikov2020}. { We also compared the performance of ROCK2 with the implicit method TR-BDF2 
for all the reaction-diffusion systems (see Table \ref{tab:performance_comparison}). While the observed number of r.h.s. evaluations was larger in all considered cases for ROCK2, the difference was not significant enough to outweigh the additional cost of linear system solves required by the TR-BDF2 method, especially in the context of batched simulation on GPUs.}

\begin{table}[t]
\centering
    \begin{tabularx}{0.7\textwidth}
    {|c||c|c||c|c|c|} 
        \cline{1-6}
          & \multicolumn{2}{|c||}{\textbf{ROCK2}} & \multicolumn{3}{c|}{\textbf{TR-BDF2}}  \\
          \cline{1-6}
         \multirow{2}{*}{} & time  & total \# r.h.s.  & time & total \# r.h.s.  & total \# linear    \\
          & steps & evaluations & steps & evaluations & system solves  \\\cline{1-6}
          \textbf{FHN} & 35 & 2430 & 60 & 1631 & 466 \\\cline{1-6}
          \textbf{GM} & 41 & 1462 & 52 & 1085 & 310 \\\cline{1-6}
          \textbf{BZ} & 38 & 1392 & 76 & 1337 & 382 \\\cline{1-6}
          \textbf{RD} & 105 & 8337 & 220 & 5411 & 1546  \\\cline{1-6}
    \end{tabularx}
    \caption{ Performance comparison between the explicit method ROCK2 and the implicit method TR-BDF2 for reaction-diffusion models. The abbreviations FHN, GM, BZ, and RD denote respectively the  FitzHugh-Nagumo model, the Gierer-Meinhardt system, the Brusselator 
system, and the reaction-diffusion-ODE model.}
    \label{tab:performance_comparison}
\end{table}

\subsubsection{Mechano-chemical models}
\label{sec:TR-BDF2}

The mechano-chemical model equations \eqref{eq:tissue}-\eqref{eq:constraint} contain the conservation law for the global arc length $S(u)$, which represents the tissue incompressibility requirement. To properly handle this constraint in  numerical simulations of the discretised equations \eqref{eq:mc:MOL}-\eqref{eq:mc:arc}, we employ implicit methods, in particular the TR-BDF2 (Trapezoidal Rule -- Backward Differentiation Formula of second order) method, which has shown good results for simulating the transient behaviour of silicon devices and circuits \cite{bank1985}. Each time step of the method $h_n$ is divided into two stages $h_{n,1} = \gamma h_n$ and $h_{n,2} = (1-\gamma) h_n$, $\gamma \in (0,1)$, which are handled with the trapezoidal and the backward differentiation rules, respectively. 


Let us define by $\bm{F}(\bm{w})$ the right-hand side of the MoL system \eqref{eq:mc:MOL} and by $S(\bm{w})$ the discretised arc length \eqref{eq:mc:arc}. Here, $\bm{w}\in\mathbb{R}^{m}$, where $m=2M_\text{dim}$ is the state vector of the system. Denote the approximation of $\bm{w}$ at time step $t_n$, $n \geq 1$ by $\bm{w}_n$. Then, the two stages of the method can be formulated as follows:
\begin{align}
	\bm{w}_{n+\gamma} &= u_n + \dfrac{\gamma h_n}{2}(F(\bm{w}_n)+F(\bm{w}_{n+\gamma})),
	\label{eq:trbdf2:stage1}\\[0.5ex]
	\alpha_2 \bm{w}_{n+1} & = \alpha_1 \bm{w}_{n+\gamma} - \alpha_0 \bm{w}_n + F(\bm{w}_{n+1}),
	\label{eq:trbdf2:stage2}
\end{align}
where 
\[
\alpha_0 = \dfrac{1-\gamma}{h_n},
\quad
\alpha_1 = \dfrac{1}{\gamma(1-\gamma)h_n},
\quad
\alpha_2 = \dfrac{2-\gamma}{(1-\gamma)h_n}.
\]

Equations \eqref{eq:trbdf2:stage1} and \eqref{eq:trbdf2:stage2} can be rewritten as two systems of non-linear equations:
\begin{equation}
\begin{array}{rl}
\bm{F}_1(\bm{w}_{n+\gamma}) &= u_n + \dfrac{\gamma h_n}{2}(F(\bm{w}_n)+F(\bm{w}_{n+\gamma})) - \bm{w}_n = 0, \\[0.5ex]
\bm{F}_2(\bm{w}_{n+1}) &= \alpha_1 \bm{w}_{n+\gamma} - \alpha_0 \bm{w}_n + F(\bm{w}_{n+1}) - \alpha_2 \bm{w}_{n+1} = 0,
\end{array}
\label{eq:trbdf2:sys}
\end{equation}
which must be solved with respect to $\bm{w}_{n+\gamma}$ and $\bm{w}_{n+1}$ during each time step, respectively. Additionally, the method must preserve the global arc length, which results in the requirement $S(\bm{w}_n) = S(\bm{w}_{n+\gamma}) = S(\bm{w}_{n+1})$. To properly take this constraint into account, instead of \eqref{eq:trbdf2:sys} we consider the extended system of equations:
\[
\begin{array}{rl}
\bm{\tilde{F}}_1(\bm{w}_{n+\gamma};\lambda_{n+\gamma}) &= (\bm{F}_1(\bm{w}_{n+\gamma};\lambda_{n+\gamma}); S(\bm{w}_{n+\gamma}) - S(\bm{w}_n)) = 0, \\[0.5ex]
\bm{\tilde{F}}_2(\bm{w}_{n+1};\lambda_{n+1}) &= (\bm{F}_2(\bm{w}_{n+1};\lambda_{n+1}); S(\bm{w}_{n+1}) - S(\bm{w}_{n+\gamma})) = 0,
\end{array}
\]
which now must be solved with respect to $\bm{\tilde{w}}_{n+\gamma} = (\bm{w}_{n+\gamma};\lambda_{n+\gamma})$ and $\bm{\tilde{w}}_{n+1}=(\bm{w}_{n+1};\lambda_{n+1})$, where $\lambda_{n+\gamma},\lambda_{n+1}\in\mathbb{R}$ are Lagrange multipliers, appearing in the right-hand side of system \eqref{eq:mc:MOL}.
In numerical simulations, we set $\gamma = 2 - \sqrt{2}$, which has the advantage that the Jacobians for both systems have the same form. As a result, we can search for $\bm{\tilde{w}}_{n+\gamma}$ and $\bm{\tilde{w}}_{n+1}$ at the same time, by using a quasi-Newton method. We define the initial guesses for the unknowns $\bm{\tilde{w}}_{n+\gamma}^{(0)} = \bm{\tilde{w}}_n$,  $\bm{\tilde{w}}_{n+1}^{(0)} = \bm{\tilde{w}}_n$, and compute iterates $\bm{\tilde{w}}_{n+\gamma}^{(k)}$ and $\bm{\tilde{w}}_{n+1}^{(k)}$, $k=1,2,\ldots$ as follows:
\[
\bm{\tilde{w}}_{n+\gamma}^{(k)} = \bm{\tilde{w}}_{n+\gamma}^{(k-1)} + \delta\bm{\tilde{w}}_{n+\gamma}^{(k)},
\quad
\bm{\tilde{w}}_{n+1}^{(k)} = \bm{\tilde{w}}_{n+1}^{(k-1)} + \delta\bm{\tilde{w}}_{n+1}^{(k)},
\]
where vectors $\delta\bm{\tilde{w}}_{n+\gamma}^{(k)}$ and $\delta\bm{\tilde{w}}_{n+1}^{(k)}$ are determined as solutions of the linear systems
\[
\bm{J} \delta\bm{\tilde{w}}_{n+\gamma}^{(k)} = - \bm{\tilde{F}}_1(\bm{\tilde{w}}_{n+\gamma}^{(k-1)}),
\quad
\bm{J} \delta\bm{\tilde{w}}_{n+1}^{(k)} = - \bm{\tilde{F}}_2(\bm{\tilde{w}}_{n+1}^{(k-1)}).
\]
Here matrix $\bm{J}$ is used as an approximation of the Jacobian $\text{\bf D} \bm{\tilde{F}}_1$ at both $\bm{\tilde{w}}_{n+\gamma}^{(k-1)}$ and $\bm{\tilde{w}}_{n+1}^{(k-1)}$. Initially, we set $J = \text{\bf D} \bm{\tilde{F}}_1(\bm{\tilde{w}}_{n+1}^{(0)})$, and update this matrix once per 50 iterations. After every update, we compute its LU-decomposition, which we later use to solve linear systems. The stopping criteria for iterates are
\[
\| \bm{\tilde{F}}_1(\bm{\tilde{w}}_{n+\gamma}^{(k)}) \|_2 < a_{tol},
\quad
\| \bm{\tilde{w}}_{n+\gamma}^{(k)} - \bm{\tilde{w}}_{n+\gamma}^{(k-1)}  \|_2 < r_{tol},
\]
and
\[
\| \bm{\tilde{F}}_1(\bm{\tilde{w}}_{n+1}^{(k)}) \|_2 < a_{tol},
\quad
\| \bm{\tilde{w}}_{n+1}^{(k)} - \bm{\tilde{w}}_{n+1}^{(k-1)}  \|_2 < r_{tol},
\]
respectively. In numerical experiments, we set $a_{tol}=10^{-3}$ and $r_{tol}=10^{-3}$.

The leading-order term of the local truncation error $\tau$  is given as follows:
\[
\tau = \dfrac{3\gamma^2 - 4\gamma + 2}{6(1-\gamma)^2}\left[ \bm{w}_{n+1} - \dfrac{1}{\gamma^2}\bm{w}_{n+\gamma} + \dfrac{1-\gamma^2}{\gamma^2}\bm{w}_n + h_n\dfrac{1-\gamma}{\gamma} \bm{F}(\bm{w}_n) \right]
\]
and the optimal step size is found using the formula 
\[
h_{opt} = h_n \sqrt[3]{\dfrac{\varepsilon_R}{\|\tau_R\|_2}},
\quad
\tau_{R,j} = \dfrac{|\tau_j|}{|w_{n+1,j}|+\delta},
\]
where in numerical simulations we set $\varepsilon_R=10^{-3}$ and $\delta=10^{-3}$.
\end{appendix}

\end{document}